\documentclass[a4paper]{amsart}
\usepackage{eucal}
\usepackage[all]{xy}
\usepackage[latin1]{inputenc}

 {\par\bigskip\centerline{{\color{red}\qquad\hrulefill\ #1 \hrulefill\qquad}}%
 \par\medskip\noindent\ignorespaces}%
 {\par\centerline{{\color{red}\qquad\hrulefill\qquad}}\par\bigskip}

\newtheorem{contar}{contar}
\newtheorem{proposition}[contar]{Proposition}

\newtheorem{theorem}[contar]{Theorem}
\newtheorem{definition}[contar]{Definition}
\newtheorem{corollary}[contar]{Corollary}
\newtheorem{remarkth}[contar]{Remark}
\newenvironment{remark}{\begin{remarkth}\upshape}{\hfill$\diamond$\end{remarkth}}

\newcommand{\qquand}{\qquad\text{and}\qquad}
\newcommand{\quand}{\quad\text{and}\quad}

\newcommand{\spV}{^{\scriptscriptstyle V}}

\newcommand{\spi}{^{\scriptscriptstyle(1)}}

\newcommand{\set}[2]{\left\{\,#1\left.\vphantom{#1#2}\,\right\vert\,#2\,
                \right\}}
\newcommand{\map}[3]{#1\colon#2\rightarrow#3}
\renewcommand{\sec}[1]{\operatorname{Sec}(#1)}                
\let\Sec\sec
\newcommand{\vectorfields}[1]{\mathfrak{X}(#1)}

\newcommand{\id}{{\operatorname{id}}}
\newcommand{\tr}{\operatorname{tr}}
\newcommand{\Lin}{\operatorname{Lin}}
\newcommand{\ext}[2][]{\bigwedge\nolimits^{#1}{#2}}
\newcommand{\cinfty}[1]{C^{\scriptscriptstyle\infty}(#1)}   
\newcommand{\pd}[2]{\frac{\partial#1}{\partial#2}}    
\let\ds\displaystyle
\newcommand{\dpd}[2]{{\ds\pd{#1}{#2}}}

\newcommand{\at}[1]{\Big\vert_{#1}}

\newcommand{\pai}[2]{\langle{#1},{#2}\rangle}         
\newcommand{\pb}{^\star}

\def\seq 0->#1-#2->#3-#4->#5->0{0\longrightarrow
            #1\maparrow{#2}#3\maparrow{#4}#5\longrightarrow 0}
\newcommand{\maparrow}[1]{\mathrel{\mathop{\longrightarrow}\limits^{#1}}}
\def\hook{\mathop{\hbox to 6pt{\hrulefill}
                      \hbox{\vrule\phantom{\vbox to 7pt{}}}}}
\newcommand{\Real}{\mathbb{R}}
\let\R\Real
\newcommand{\bs}[1]{\boldsymbol{#1}}

\makeatletter

\newcommand\prol{\@ifstar{\@proldf}{\@prolpf}}  

\def\@prolpf{\@ifnextchar[{\@prolpf@wrt}{\@prolpf@}}
\def\@prolpf@wrt[#1]#2{\@ifnextchar[{\@prolpf@wrt@at{#1}{#2}}{\@prolpf@wrt@{#1}{#2}}}
\def\@prolpf@wrt@at#1#2[#3]{\prolsymbol^{#1}_{#3}#2}
\def\@prolpf@wrt@#1#2{\prolsymbol^{#1}#2}
\def\@prolpf@#1{\@ifnextchar[{\@prolpf@at{#1}}{\@prolpf@@{#1}}}
\def\@prolpf@at#1[#2]{\prolsymbol_{#2}#1}
\def\@prolpf@@#1{\prolsymbol#1}

\def\@proldf{\@ifnextchar[{\@proldf@wrt}{\@proldf@}}
\def\@proldf@wrt[#1]#2{\@ifnextchar[{\@proldf@wrt@at{#1}{#2}}{\@proldf@wrt@{#1}{#2}}}
\def\@proldf@wrt@at#1#2[#3]{\prolsymbol^{*#1}_{#3}#2}
\def\@proldf@wrt@#1#2{\prolsymbol^{*#1}#2}
\def\@proldf@#1{\@ifnextchar[{\@proldf@at{#1}}{\@proldf@@{#1}}}
\def\@proldf@at#1[#2]{\prolsymbol^*_{#2}#1}
\def\@proldf@@#1{\prolsymbol^*#1}

\def\prolsymbol{\CMcal{T}}

\makeatother

\let\ol\overline
\let\ul\underline
\newcommand{\prEM}{\tau^E_M}
\newcommand{\prFN}{\tau^F_N}
\newcommand{\prEF}{\ol{\pi}}
\newcommand{\prMN}{\ul{\pi}}
\newcommand{\Lpi}[1][]{\mathcal{L}_{#1}\pi}
\newcommand{\Jpi}[1][]{\mathcal{J}_{#1}\pi}
\newcommand{\Vpi}[1][]{\mathcal{V}_{#1}\pi}
\newcommand{\JJpi}[1][]{\mathcal{J}_{#1}\pi_1}
\newcommand{\Jtpi}[1][]{\mathcal{J}^2_{#1}\pi}
\newcommand{\prLM}{\ul{\tilde{\pi}_{10}}}
\newcommand{\prVM}{\ul{\bs{\pi}_{10}}}
\newcommand{\prJM}{\ul{\pi_{10}}}
\newcommand{\prJN}{\ul{\pi_1}}
\newcommand{\prTJF}{\ol{\pi_1}}

\newcommand{\dLpi}[1][]{\mathcal{L}_{#1}^*\!\pi}
\newcommand{\dJpi}[1][]{\mathcal{J}_{#1}^\dag\pi}
\newcommand{\dVpi}[1][]{\mathcal{V}_{#1}^*\!\pi}
\newcommand{\prdJM}{\prJM^\dag}
\newcommand{\prdVM}{\prVM^*}

\newcommand{\prTdJE}{\ol{\pi_{10}}{}^\dag}
\newcommand{\pr}{\mathrm{pr}}

\newcommand{\Mor}[1]{\mathcal{M}(#1)} 
\newcommand{\Adm}[1]{\operatorname{Sec}^{\operatorname{adm}}(#1)}
\newcommand{\X}{\mathcal{X}}   
\newcommand{\V}{\mathcal{V}}   
\renewcommand{\P}{\mathcal{P}}

\def\p#1{\mathaccent19{#1}}

\renewcommand{\v}[1][]{{v^{#1}}}

\renewcommand{\L}{\mathcal{L}} 
\newcommand{\leg}{\widehat{\mathcal{F}}_\L}
\newcommand{\legd}{\mathcal{F}_\L}

\def\emph#1{{\bfseries\itshape{#1}}}

\begin{document}

\title[Classical Field theory on Lie algebroids]{Classical Field theory on Lie algebroids: multisymplectic formalism}

\author{Eduardo Mart\'{\i}nez}
\address{Eduardo Mart\'{\i}nez:
Departamento de Matem\'atica Aplicada,
Facultad de Ciencias,
Universidad de Zaragoza,
50009 Zaragoza, Spain}
\email{emf@unizar.es}

\thanks{Partial financial support from MICYT grant  BFM2003-02532 is acknowledged}

\keywords{Lie algebroids, jet bundles, Lagrangian and Hamiltonian Field theory, multisymplectic forms}
\subjclass[2000]{58A20, 70S05, 49S05, 35F20, 53C15, 53D99, 58H99}

\date{October 15, 2004 (Preliminary versions: June 2001, December 2002)}

\begin{abstract}
The jet formalism for Classical Field theories is extended to the setting of Lie algebroids. We define the analog of the concept of jet of a section of a bundle and we study some of the geometric structures of the jet manifold. When a Lagrangian function is given, we find the equations of motion in terms of a Cartan form canonically associated to the Lagrangian. The Hamiltonian formalism is also extended to this setting and we find the relation between the solutions of both formalism. When the first Lie algebroid is a tangent bundle we give a variational description of the equations of motion. In addition to the standard case, our formalism includes as particular examples the case of systems with symmetry (covariant Euler-Poincaré and Lagrange Poincaré cases), variational problems for holomorphic maps, Sigma models or Chern-Simons theories. One of the advantages of our theory is that it is based in the existence of a multisymplectic form on a Lie algebroid.
\end{abstract}

\maketitle


\section{Introduction}

The idea of using Lie algebroids in Mechanics is due to Weinstein~\cite{Weinstein}. He showed that it is possible to give a common description of the most interesting classical mechanical systems by using the geometry of Lie algebroids (and their discrete analogs using Lie groupoids). The theory includes as particular cases systems with holonomic constraints, systems defined on Lie algebras, systems with symmetry and systems defined on semidirect products, described by the Euler-Lagrange equations (also known in each case as constrained Euler-Lagrange, Euler-Poincaré, Lagrange-Poincaré or Hammel equations~\cite{CeMaRa}). 

Following the ideas by Klein~\cite{Klein}, a geometric formalism  was introduced in~\cite{LMLA}, showing that the Euler-Lagrange equations are obtained be means of a symplectic equation. Here symplectic is to be understood in the context of Lie algebroids, i.e. a regular skew-symmetric bilinear form which is closed with respect to the exterior differential operator defined by the Lie algebroid structure. Also a Hamiltonian symplectic description was defined in~\cite{HMLA} (see also~\cite{Medina,LeMaMa}).

Later~\cite{SaMeMa1,MaMeSa} this theory was extended to time-dependent Classical Mechanics  by introducing a generalization of the notion of Lie algebroid when the bundle is no longer a vector bundle but an affine bundle. The particular case considered in~\cite{SaMeMa2} is the analog of a first jet bundle of a bundle $M\to\R$, and hence it is but a particular case of a first order Field Theory for a 1-dimensional space-time. Therefore, it is natural to investigate whether it is possible to extend our formalism to the case of a general Field Theory, where the space-time manifold is no longer one dimensional, and where the field equations are defined in terms of a variational problem. 

The today standard formalism in classical Field Theory is the multisymplectic approach. It was implicit in the work of De Donder~\cite{DeDonder} and Lepage~\cite{Lepage}, and it was rediscovered later on in the context of relativistic field theories~\cite{GoSt,Garcia,KiTu,BiSnFi}. See also~\cite{CaCrIb,GIMMSY1, GIMMSY2,Gotay,Barna-L,Barna-H,Barna-MC1,Barna-MC2,HeKo,HeKo2,PaRo}. A recent renewed interest in multisymplectic geometry is in part motivated by the discovery of numerical integrators which preserve the multisymplectic form~\cite{BrRe,MaPaSh,MaPeShWe}

There exists other approaches in the literature, as it is the case of the polysymplectic formalism in its various versions: Günther's formalism~\cite{Gunther,LeMeSa, LeMeOuRoSa,MuReSa} and that of Sardanashvily and coworkers~\cite{MaSa, GiMaSa, Sardanashvily}. See also~\cite{Kanatchikov,Norris}. They will not be considered in this paper.

\medskip

Generally speaking, there are three different but closely related aspects in the analysis of first order field theories:
\begin{itemize}
\item\textsl{variational approach}, which leads to the Euler-Lagrange equations,
\item\textsl{multisymplectic formalism} on the first jet bundle, and
\item\textsl{infinite-dimensional approach} on the space of Cauchy data.
\end{itemize}

In a previous paper~\cite{DGA04} we have studied several aspects of the variational approach to the theory formulated in the context of Lie algebroids, a review of which is presented here in a more intrinsic way. The aim of this paper is to develop the multisymplectic formalism for those variational problems defined over Lie algebroids, and we will leave for the future the study of the infinite-dimensional formalism.

\medskip

The standard geometrical approach to the Lagrangian description of first order classical field theories~\cite{CaCrIb,Barna-L,GIMMSY1} is based on the canonical structures on the first order jet bundle~\cite{Saunders} of a fiber bundle whose sections are the fields of the theory. Thinking of a Lie algebroid as a substitute for the tangent bundle to a manifold, the analog of the field bundle to be considered here is a surjective morphism of Lie algebroids $\map{\pi}{E}{F}$. 

In the standard theory, a 1-jet of a section of a bundle is but the tangent map to that section at the given point, and therefore it is a linear map between tangent spaces which has to be a section of the tangent of the projection map. In our theory, the analog object is a linear map from a fiber of $F$ to a fiber of $E$ which is a section of the projection $\pi$. The space $\Jpi$ of these maps has the structure of an affine bundle and carries similar canonical structures to those of the first jet bundle, and it is in that space where our theory is based. In particular, in the affine dual of $\Jpi$ there exists a canonical multisymplectic form, which allows to define a Hamiltonian formalism. 

One relevant feature of our theory is that it is a multisymplectic theory,  i.e. in both the (regular) Lagrangian and the Hamiltonian approaches, the field equations are expressed in term of a multisymplectic form $\Omega$, as the equations for some sections $\Psi$ satisfying an equation of the form $\Psi\pb i_X\Omega=0$. In this way it is possible to extend the category of multisymplectic manifolds to the category of multisymplectic Lie algebroids, where it is possible to study reduction within this category. 

The situation is to be compared with the so called covariant reduction theory~\cite{CaRaSh,CaGaRa,CaRa}.  There the field equations are obtained by reducing the variational problem, i.e. by restricting the variations to those coming from variations for the original unreduced problem. Therefore, that equations are of different form in each case, Euler-Poincaré, Lagrange-Poincaré for a system with symmetry and (a somehow different) Lagrange-Poincaré for semidirect products. In contrast, our theory includes all these cases as particular cases and the equations are always of the same form. This fact must be important in the general theory of reduction.

\medskip

The organization of the paper is as follows. In section~\ref{preliminaries} we will recall some differential geometric structures related to Lie algebroids, as the exterior differential, the tangent prolongation of a bundle with respect to a Lie algebroid, the tangent prolongation of a bundle map with respect to a morphism of Lie algebroids, and the flow defined by a section. In section~\ref{jetoids} we define the manifold of jets in this generalized sense, and we will show that it is (the total space of) an affine bundle. The anchor and the bracket of the Lie algebroid introduce more structure which is presented there. In section~\ref{2-jets} we consider repeated jets. This is necessary since the equations in our theory are, as in the standard case, second order partial differential equations (generally speaking). This leads us to the definition of holonomic and semiholonomic jets. The relation of jets on Lie algebroids and ordinary jets is studied in section~\ref{morphisms} where we also show how to associate a partial differential equation to a submanifold of $\Jpi$. We pay special attention to the conditions relating holonomic jets and morphisms of Lie algebroids. In section~\ref{Lagrangian} we introduce the Lagrangian description of a Field Theory defined by a Lagrangian on $\Jpi$. The equations of motion are expressed in terms of a Cartan form and of its differential, the multisymplectic form. In section~\ref{variational} we show that, in the case $F=TN$, the equations obtained by the multisymplectic formalism coincide with the Euler-Lagrange equations determined by constrained variational calculus. The Hamiltonian counterpart is studied in section~\ref{Hamiltonian}. We show the relation between the solutions in the Hamiltonian and in the Lagrangian approaches. In section~\ref{examples} we show some examples where the theory can be applied and to end the paper we discuss some open problems and future work.

We finally mention that related ideas to those explained in this paper have been recently considered by Strobl and coworkers~\cite{Strobl,BoKoSt} in the so called off-shell theory. In this respect our theory should be considered as the on-shell counterpart for the cases considered there.

\subsection*{Notation}
All manifolds and maps are taken in the $C^\infty$ category. The set of smooth functions on a manifold $M$ will be denoted by $\cinfty{M}$. The set of smooth vectorfields on a manifold $M$ will be denoted $\vectorfields{M}$. The set of smooth sections of a fiber bundle $\map{\pi}{B}{M}$  will be denoted $\sec{\pi}$ or $\sec{B}$ when there is no possible confusion. 

The tensor product of a vector bundle $\map{\tau}{E}{M}$ by itself $p$ times will be denoted by $\map{\tau^{\otimes p}}{E^{\otimes p}}{M}$, and similarly the exterior power will be denoted by  $\map{\tau^{\wedge p}}{E^{\wedge p}}{M}$. The set of sections of the dual $\map{\tau^*{}^{\wedge p}}{(E^*)^{\wedge p}}{M}$ will be denoted by $\ext[p]{E}$. For $p=0$ we have $\ext[0]{E}=\cinfty{M}$.

By a bundle map from $\map{\pi}{P}{M}$ to $\map{\pi'}{P'}{M'}$ we mean a pair $\Phi=(\ol{\Phi},\ul{\Phi})$ with $\map{\ol{\Phi}}{P}{P'}$ and $\map{\ul{\Phi}}{M}{M'}$ satisfying $\pi'\circ\ol{\Phi}=\ul{\Phi}\circ\pi$. The map $\ul{\Phi}$ is said to be the base map. If $P$ and $P'$ are vector bundles and $\ol{\Phi}$ is fiberwise linear we say that $\Phi$ is a vector bundle map.

Let $\map{\tau}{E}{M}$ and $\map{\tau'}{E'}{M'}$ be two vector bundles. Consider two points $m\in M$ and $m'\in M'$, and a linear map $\map{\phi}{E_m}{E'_{m'}}$. The transpose or dual map is a map $\map{\phi\pb}{E'{}^*_{m'}}{E^*_m}$ given by  $\pai{\phi\pb\alpha}{a}=\pai{\alpha}{\phi(a)}$, for every $\alpha\in E'{}^*_{m'}$. This action extends to covariant tensors as usual. If moreover $\phi$ is invertible, then we can extend $\phi\pb$ also to contravariant tensors. We just need to define it for vectors: if $\zeta\in E'_{m'}$ then $\phi\pb\zeta=\phi^{-1}(\zeta)\in E_m$.

The notion of pullback of a tensor field by a vector bundle map needs some attention.  Given a vector bundle map $\Phi=(\ol{\Phi},\ul{\Phi})$ from the vector bundle $\map{\tau}{E}{M}$ to the vector bundle $\map{\tau'}{E'}{M'}$ we define the pullback operator as follows. For every section $\beta$ of the $p$-covariant tensor bundle $\otimes ^pE'{}^*\to M'$ we define the section $\Phi\pb\beta$ of $\otimes^p E^*\to M$ by
$$
(\Phi\pb\beta)_m(a_1,a_2,\ldots,a_p)=
\beta_{\ul{\Phi}(m)}\bigl(\ol{\Phi}(a_1),\ol{\Phi}(a_2),\ldots,\ol{\Phi}(a_p)\bigr),
$$
for every $a_1,\ldots,a_p\in E_m$. The tensor $\Phi\pb\beta$ is said to be the pullback of $\beta$ by $\Phi$. For a function $f\in\cinfty{M'}$ (i.e. for $p=0$) we just set $\Phi\pb f=\ul{\Phi}^*f=f\circ\ul{\Phi}$. 

It follows that $\Phi\pb(\alpha\otimes\beta)=(\Phi\pb\alpha)\otimes(\Phi\pb\beta)$. When $\Phi$ is fiberwise invertible, we define the pullback of a section $\sigma$ of $E'$ by $(\Phi\pb\sigma)_m=(\ol\Phi_m)^{-1}(\sigma(\ul\Phi(m))$, where $\ol\Phi_m$ is the restriction of $\ol\Phi$ to the fiber $E_m$ over $m\in M$. Thus, the pullback of any tensor field is defined.

In the case of the tangent bundles $E=TM$ and $E'=TM'$, and the tangent map $\map{T\varphi}{TM}{TM'}$ of a map $\map{\varphi}{M}{M'}$ we have that $(T\varphi,\varphi)\pb\beta=\varphi^*\beta$ is the ordinary pullback by $\varphi$ of the covariant tensor field $\beta$ on $M'$. Notice the difference between the symbols ${}\pb$ (star) and ${}^*$ (asterisque). 


\section{Preliminaries on Lie algebroids}
\label{preliminaries}

In this section we introduce some well known notions and develop new
concepts concerning Lie algebroids that are necessary for the further
developments. We refer the reader to~\cite{CaWe,Mackenzie,Mackenzie2,HiMa} for
details about Lie groupoids, Lie algebroids and their role in
differential geometry. While our main interest is in Lie algebroids, we will consider as an intermediate step the category of anchored vector bundles~\cite{CaLa,PoPo}.

\subsection*{Lie algebroids}
Let $M$ be an $n$-dimensional manifold and let $\map{\tau}{E}{M}$ be a vector bundle. A vector bundle map $\map{\rho}{E}{TM}$ over the identity is called an anchor map. The vector bundle $E$ together with an anchor map $\rho$ is said to be an \emph{anchored vector bundle}. A structure of \emph{Lie algebroid} on $E$ is given by a Lie algebra structure on the $\cinfty{M}$-module of sections of the bundle, $(\Sec{E},[\cdot\,,\cdot])$, together with an anchor map, satisfying the compatibility condition
$$
[\sigma,f\eta] = f[\sigma,\eta] + \bigl(
\rho(\sigma)f \bigr) \eta .
$$
Here $f$ is a smooth function on $M$, $\sigma$, $\eta$ are
sections of $E$ and we have denoted by $\rho(\sigma)$ the vector field
on $M$ given by $\rho(\sigma)(m)=\rho(\sigma(m))$. From the compatibility
condition and the Jacobi identity, it follows that the map
$\sigma\mapsto\rho(\sigma)$ is a Lie algebra homomorphism from
$\Sec{E}$ to $\mathfrak{X}(M)$. 

It is convenient to think of a Lie algebroid, and more generally an anchored vector bundle, as a substitute of the tangent bundle of $M$. In this way, one regards an element $a$ of $E$ as a generalized velocity, and the actual velocity $v$ is obtained when applying the anchor to $a$, i.e., $v=\rho(a)$.
A curve $\map{a}{[t_0,t_1]}{E}$ is said to be \emph{admissible} if $\dot{m}(t)=\rho(a(t))$, where $m(t)=\tau(a(t))$ is the base curve.

When $E$ carries a Lie algebroid structure, the image of the anchor map, $\rho(E)$, defines an integrable smooth generalized distribution on $M$.  Therefore, $M$ is foliated by the integral leaves of $\rho(E)$, which are called the leaves of the Lie algebroid. It follows that $a(t)$ is admissible if and only if the curve $m(t)$ lies on a leaf of the Lie algebroid, and that two points are in the same leaf if and only if they are connected by (the base curve of) an admissible curve.

A Lie algebroid is said to be transitive if it has only one leaf, which is obviously equal to $M$. It is easy to see that $E$ is transitive if and only if $\rho$ is surjective. If $E$ is not transitive, then the restriction of the Lie algebroid to a leaf $L\subset M$, $E_{|L}\to L$ is transitive.

Given local coordinates $(x^i)$ in the base manifold $M$ and a local basis $\{e_\alpha\}$ of sections of $E$, we have local coordinates $(x^i,y^\alpha)$ in $E$. If $a\in E_m$ is an element then we can write $a=a^\alpha e_\alpha(m)$ and thus the coordinates of $a$ are $(m^i,a^\alpha)$, where $m^i$ are the coordinates of the point $m$. The anchor map is locally determined by the local functions $\rho^i_\alpha$ on $M$ defined by $\rho(e_\alpha)=\rho^i_\alpha(\partial/\partial x^i)$. In addition, for a Lie algebroid, the Lie bracket is determined by the functions  $C^\alpha_{\beta\gamma}$ defined by $[e_\alpha,e_\beta]=C^\gamma_{\alpha\beta}e_\gamma$. The functions $\rho^i_\alpha$ and $C^\alpha_{\beta\gamma}$ are said to be the structure functions of the Lie algebroid in this coordinate system. They satisfy the following relations
\begin{align*}
  \rho^j_\alpha\pd{\rho^i_\beta}{x^j} -
  \rho^j_\beta\pd{\rho^i_\alpha}{x^j} = \rho^i_\gamma
  C^\gamma_{\alpha\beta},
  \quand
  \sum_{\mathrm{cyclic}(\alpha,\beta,\gamma)} \left[\rho^i_\alpha\pd{
      C^\nu_{\beta\gamma}}{x^i} + C^\mu_{\alpha\nu}
    C^\nu_{\beta\gamma}\right]=0,
\end{align*}
which are called the structure equations of the Lie algebroid. 

\subsection*{Exterior differential}
The anchor $\rho$ allows to define the differential of a function on the base manifold with respect to an element $a\in E$. It is given by
$$
df(a)=\rho(a)f.
$$
It follows that the differential of $f$ at the point $m\in M$ is an element of $E_m^*$. 

Moreover, a structure of Lie algebroid on $E$ allows to extend the differential to sections of the bundle $\ext[p]{E}$ which we will call just $p$-forms. If $\omega\in\sec{\ext[p]{E}}$ then $d\omega\in\sec{\ext[p+1]{E}}$ is defined by 
\begin{align*}
d\omega(\sigma_0,\sigma_1,\ldots,\sigma_p)
&=
\sum_i\rho(\sigma_i)(-1)^i
  \omega(\sigma_0,\ldots,\widehat{\sigma_i},\ldots,\sigma_p)+\\
&\qquad{}+
\sum_{i<j}(-1)^{i+j}
\omega([\sigma_i,\sigma_j],\sigma_0,\ldots,\widehat{\sigma_i},\ldots,\widehat{\sigma_j},\ldots,\sigma_p).
\end{align*}
It follows that $d$ is a cohomology operator, that is $d^2=0$.

Locally the exterior differential is determined by
$$
dx^i=\rho^i_\alpha e^\alpha
\qquand
de^\gamma=-\frac{1}{2}C^\gamma_{\alpha\beta}e^\alpha\wedge e^\beta.
$$
In the rest of the paper the symbol $d$ will refer to the exterior differential on the Lie algebroid $E$ and not to the ordinary exterior differential on a manifold.

The usual Cartan calculus extends to the case of Lie algebroids (See~\cite{Nijenhuis}). For every section $\sigma$ of $E$ we have a derivation $i_\sigma$ (contraction) of degree $-1$ and a derivation $d_\sigma=i_\sigma\circ d+d\circ i_\sigma$ (Lie derivative) of degree $0$. Since $d^2=0$ we have that $d_\sigma\circ d=d\circ d_\sigma$. Nevertheless, notice that there are derivations of degree $0$ that commutes with the exterior differential $d$ but are not Lie derivatives.

\subsection*{Admissible maps and morphisms}
Let $\map{\tau}{E}{M}$ and $\map{\tau'}{E'}{M'}$ be two
anchored vector bundles, with anchor maps $\map{\rho}{E}{TM}$ and
$\map{\rho'}{E'}{TM'}$.  Let $\Phi=(\ol{\Phi},\ul{\Phi})$ be a vector bundle map, that is
$\map{\ol{\Phi}}{E}{E'}$ is a fiberwise linear map over $\map{\ul{\Phi}}{M}{M'}$. We will say that $\Phi$ is admissible if it maps admissible curves into admissible curves. It follows that $\Phi$ is admissible if and only if 
$T\ul{\Phi}\circ\rho = \rho'\circ\ol{\Phi}$. This condition can be conveniently expressed in terms of the exterior differential as follows. The map $\Phi$ is \emph{admissible} if and only if $\Phi\pb df=d\Phi\pb f$ for every function $f\in\cinfty{M'}$. 

If $E$ and $E'$ are Lie algebroids, then we say that $\Phi$ is a \emph{morphism} if $\Phi\pb d\theta=d\Phi\pb\theta$ for every $\theta\in\Sec{\ext{E}}$. Obviously, a morphism is an admissible map. It is clear that an admissible map can be properly called a morphism of anchored vector bundles; nevertheless since our main interest is in Lie algebroids we will reserve the word \textsl{morphism} for a morphism of Lie algebroids.

\medskip

Let $(x^i)$ be a local coordinate system on $M$ and $(x'{}^i)$ a local coordinate system on $M'$. Let $\{e_\alpha\}$ and $\{e'_\alpha\}$ be local basis of sections of $E$ and $E'$, respectively, and $\{e^\alpha\}$ and $\{e'{}^\alpha\}$ the dual basis. The bundle map $\Phi$ is determined by the relations $\Phi\pb x'{}^i = \phi^i(x)$ and $\Phi\pb e'{}^\alpha = \phi^\alpha_\beta e^\beta$ for certain local functions $\phi^i$ and $\phi^\alpha_\beta$ on $M$. Then $\Phi$ is admissible if and only if
$$
  \rho^j_\alpha\pd{\phi^i}{x^j}=\rho'{}^i_\beta\phi^\beta_\alpha.
$$
The map $\Phi$ is a morphism of Lie algebroids if and only if 
$$
  \phi^\beta_\gamma C^\gamma_{\alpha\delta} =
    \left(\rho^i_\alpha\pd{\phi^\beta_\delta}{x^i} -
    \rho^i_\delta\pd{\phi^\beta_\alpha}{x^i}\right) +
  C'{}^\beta_{\theta\sigma}\phi^\theta_\alpha\phi^\sigma_\delta,
$$
in addition to the admissibility condition above. In these expressions $\rho^i_\alpha$, $C^\alpha_{\beta\gamma}$ are the structure functions on $E$ and $\rho'{}^i_\alpha$, $C'{}^\alpha_{\beta\gamma}$ are the structure functions on $E'$. 

\subsection*{The prolongation of a Lie algebroid} (See~\cite{LMLA,MaMeSa,PoPo,LeMaMa}) Let $\map{\pi}{P}{M}$ be a bundle with base manifold $M$. Thinking of $E$ as a substitute of the tangent bundle to $M$, the tangent bundle to $P$ is not the appropriate space to describe dynamical systems on $P$. This is clear if we note that the projection to $M$ of a vector tangent to $P$ is a vector tangent to $M$, and what one would like instead is an element of $E$, the `new' tangent bundle of $M$.

A space which takes into account this restriction is the $E$-tangent
bundle of $P$, also called the prolongation of $P$ with respect to $E$, which we denote here by $\prol[E]{P}$. It is defined as the vector
bundle $\map{\tau^E_P}{\prol[E]{P}}{P}$ whose fiber at a point $p\in P_m$
is the vector space
$$
\prol[E]{P}[p] =\set{(b,v)\in E_m\times T_pP}{\rho(b)=T_p\pi(v)}.
$$
We will frequently use the redundant notation $(p,b,v)$ to denote the element $(b,v)\in\prol{E}[p]$. In this way, the map $\tau^E_P$ is just the projection onto the first factor. The anchor of $\prol[E]{P}$ is the projection onto the third factor, that is, the map $\map{\rho^1}{\prol[E]{P}}{TP}$ given by $\rho^1(p,b,v)=v$. 

We also consider the map $\map{\prol{\pi}}{\prol[E]{P}}{E}$ defined by $\prol{\pi}(p,b,v)=b$. The bundle map $(\prol{\pi},\pi)$ will also be denoted by $\prol{\pi}$. The origin of the notation $\prol{\pi}$ will become clear later on.
An element $Z=(p,b,v) \in \prol[E]{P}$ is said to be vertical if it projects to zero, that is $\prol{\pi}(Z)=0$. Therefore it is of the form $(p,0,v)$, with $v$ a vertical vector tangent to $P$ at $p$. 

Given local coordinates $(x^i,u^A)$ on $P$ and a local basis
$\{e_\alpha\}$ of sections of $E$, we can define a local basis
$\{\X_\alpha,\V_A\}$ of sections of $\prol[E]{P}$ by
$$
\X_\alpha(p) =\Bigl(p,e_\alpha(\pi(p)),\rho^i_\alpha\pd{}{x^i}\at{p}\Bigr)
  \qquand
\V_A(p) = \Bigl(p,0,\pd{}{u^A}\at{p}\Bigr).
$$
If $(p,b,v)$ is an element of $\prol[E]{P}$, with $b=z^\alpha e_\alpha$, then $v$ is of the form $v=\rho^i_\alpha z^\alpha\pd{}{x^i}+v^A\pd{}{u^A}$, and we can write
$$
(p,b,v)=z^\alpha\X_\alpha(p)+v^A\V_A(p).
$$
Vertical elements are linear combinations of $\{\V_A\}$. 

The anchor map $\rho^1$ applied to a section $Z$  of $\prol[E]{P}$ with local expression $Z =
Z^\alpha\X_\alpha+V^A\V_A$ is the  vectorfield on $E$ whose coordinate expression is
$$
\rho^1(Z) = \rho^i_\alpha Z^\alpha \pd{}{x^i} + V^A\pd{}{u^A}.
$$

If $E$ carries a Lie algebroid structure, then so does $\prol[E]{P}$. The  Lie bracket associated can be easily defined in terms of projectable sections as follows. A section $Z$ of $\prol[E]{P}$ is said to be projectable if there exists a section $\sigma$ of $E$ such that $\prol{\pi}\circ Z=\sigma\circ\pi$.  Equivalently, a section $Z$ is projectable if and only if it is of the form $Z(p)=(p,\sigma(\pi (p)),X(p))$, for some section $\sigma$ of $E$ and some vector field $X$ on $E$. The Lie bracket of two projectable sections $Z_1$ and $Z_2$ is then given by
$$
[Z_1,Z_2](p)=(p,[\sigma_1,\sigma_2](m),[X_1,X_2](p)), \qquad p \in P,\, m=\pi(p).
$$
It is easy to see that $[Z_1,Z_2](p)$ is an element of $\prol[E]{P}[p]$ for every $p\in P$. Since any section of $\prol[E]{P}$ can be locally written as a linear combination of projectable sections, the definition of the Lie bracket for arbitrary sections of $\prol[E]{P}$ follows.

The Lie brackets of the elements of the basis are
$$
[\X_\alpha,\X_\beta]= C^\gamma_{\alpha\beta}\:\X_\gamma,
\qquad [\X_\alpha,\V_B]=0 \qquand
[\V_A,\V_B]=0.
$$
and the exterior differential is determined by
\begin{align*}
&dx^i=\rho^i_\alpha \X^\alpha
&&du^A=\V^A\\
\intertext{and}
&d\X^\gamma=-\frac{1}{2}C^\gamma_{\alpha\beta}\X^\alpha\wedge\X^\beta
&&d\V^A=0,
\end{align*}
where $\{\X^\alpha,\V^A\}$ is the dual basis to $\{\X_\alpha,\V_A\}$.

\subsection*{Prolongation of a map}
The tangent map extends to a map between the prolongations whenever we have an admissible map. Let $\Psi=(\ol{\Psi},\ul{\Psi})$ be a bundle map from the bundle $\map{\pi}{P}{M}$ to the bundle $\map{\pi'}{P'}{M'}$. Let $\Phi=(\ol{\Phi},\ul{\Phi})$ be an admissible map from $\map{\tau}{E}{M}$ to $\map{\tau'}{E'}{M'}$. Moreover, assume that both base maps coincide, $\ul{\Psi}=\ul{\Phi}$. The prolongation of $\ol{\Psi}$ with respect to $\ol{\Phi}$ is the mapping $\map{\prol[\ol{\Phi}]{\ol{\Psi}}}{\prol[E]{P}}{\prol[E']{P'}}$ defined by $$
\prol[\ol{\Phi}]{\ol{\Psi}}(p,b,v)
=(\ol{\Psi}(p),\ol{\Phi}(b),T_p\ol{\Psi}(v)).
$$
It is clear from the definition that $\prol[\Phi]{\Psi}\equiv(\prol[\ol{\Phi}]{\ol{\Psi}},\ol{\Psi})$ is a vector bundle map from $\map{\tau^E_P}{\prol[E]{P}}{P}$ to $\map{\tau^{E'}_{P'}}{\prol[E']{P'}}{P'}$. 

\begin{proposition}
The map $\prol[\Phi]{\Psi}$ is an admissible map. Moreover, $\prol[\Phi]{\Psi}$ is a morphism of Lie algebroids if and only if $\Phi$ is a morphism of Lie algebroids.
\end{proposition}
\begin{proof}
It is clearly admissible because $\rho^1\circ\prol[\ol{\Phi}]{\ol{\Psi}}=T\ol{\Psi}\circ\rho^1$. Therefore we have $d(\prol[\Phi]{\Psi})\pb F=(\prol[\Phi]{\Psi})\pb dF$ for every function $F\in\cinfty{P'}$. In order to prove that $d(\prol[\Phi]{\Psi})\pb\beta=(\prol[\Phi]{\Psi})\pb d\beta$ for every $\beta\in\sec{E'{}^*}$, we take into account that every 1-form on $\prol[E']{P'}$ can be expressed as a linear combination of exact forms $dF$ and basic forms $(\prol{\pi'})\pb\theta$, for $F\in\cinfty{P'}$ and $\theta\in\sec{E'{}^*}$. (This is clear in coordinates where $\X^\alpha=(\prol{\pi})\pb e^\alpha$ and $\V^A=dy^A$). Therefore it is enough to prove that relation for this kind of forms. 

For an exact 1-form $\beta=dF$ we have
$$
d(\prol[\Phi]{\Psi})\pb dF-(\prol[\Phi]{\Psi})\pb d dF=
dd(\prol[\Phi]{\Psi})\pb F-0=0
$$
where we have used that $\prol[\Phi]{\Psi}$ is admissible.

For a basic form $\beta= (\prol{\pi'})\pb\theta$ we have
$$
d(\prol[\Phi]{\Psi})\pb(\prol{\pi'})\pb\theta
=d(\prol{\pi'}\circ\prol[\Phi]{\Psi})\pb\theta
=d(\Phi\circ\prol{\pi})\pb\theta
=d(\prol{\pi})\pb\Phi\pb\theta
=(\prol{\pi})\pb d\Phi\pb\theta,
$$
where we have used that $\prol{\pi}$ is a morphism. On the other hand
$$
(\prol[\Phi]{\Psi})\pb d(\prol{\pi'})\pb\theta
=(\prol[\Phi]{\Psi})\pb(\prol{\pi'})\pb d\theta
=(\prol{\pi'}\circ\prol[\Phi]{\Psi})\pb d\theta
=(\Phi\circ\prol{\pi})\pb d\theta
=(\prol{\pi})\pb \Phi\pb d\theta,
$$
where we have used that $\prol{\pi'}$ is a morphism. Therefore, $\prol[\Phi]{\Psi}$ is a morphism if and only if
$$
(\prol{\pi})\pb(d\Phi\pb\theta-\Phi\pb d\theta)=0.
$$
Since $\prol{\pi}$ is surjective the above equation holds if and only if $\Phi$ is a morphism.
\end{proof}

Given local coordinates $x^i$ on $M$ and $x'{}^i$ on $M'$, local adapted coordinates $(x^i,u^A)$ on $P$ and $(x'{}^i,u'{}^A)$ on $P'$ and local basis of section $\{e_\alpha\}$ of $E$ and
$\{e'_\alpha\}$ of $E'$, the maps $\Phi$ and $\Psi$ are determined by
$\Phi\pb e'{}^\alpha=\Phi^\alpha_\beta e^\beta$ and  $\ol{\Psi}(x,u)=(\phi^i(x),\psi^A(x,u))$. Then the action of $\prol[\Phi]{\Psi}$ is given by
\begin{align*}
(\prol[\Phi]{\Psi})\pb\X'{}^\alpha
&= \Phi_\beta^\alpha\X^\beta\\
(\prol[\Phi]{\Psi})\pb\V'{}^A
&=\rho^i_\alpha\pd{\psi^A}{x^i}\X^\alpha+\pd{\Psi^A}{u^B}\V^B,
\end{align*}
Equivalently, for every $p\in P$ we have
\begin{align*}
\prol[\ol{\Phi}]{\ol{\Psi}}[p](\X_\alpha(p)) 
&= \Phi^\beta_\alpha(m)\X'_\beta(p') 
       +\rho^i_\alpha(m)\pd{\Psi^B}{x^i}(p)\V'_B(p'),  \\
\prol[\ol{\Phi}]{\ol{\Psi}}[p](\V_A(p)) 
&=\pd{\Psi^B}{u^A}(p)\,\V'_B(p'),
\end{align*}
where $m=\pi(p)$ and $p'=\Psi(p)$. 

\begin{remark}
The following special situations are relevant:
\begin{itemize}
\item  If $E=E'$ and the base map of $\Psi$ is the identity in $M$, $\ul{\Psi}=\id_M$, we can always take as $\Phi$ the identity in $E$, $\Phi=\id=(\id_E,\id_M)$. In this case, we will simplify the notation and we will write $\prol{\Psi}$ instead of $\prol[\id]{\Psi}$.

\item We can consider as the bundle $P$ the degenerate case of  $\map{\id_M}{M}{M}$. Then we have a canonical identification $E\equiv\prol[E]{M}$, given explicitly by $a\in E\equiv (m,a,\rho(a))\in\prol[E]{M}$, where $m=\tau(a)$.

\item With this identification, the map $\map{\prol{\pi}}{\prol[E]{P}}{\prol[E]{M}\equiv E}$ is just $(p,a,v)\mapsto (m,a,\rho(a))\equiv a$, which  justifies the notation $\prol{\pi}$ for the projection onto the second factor.

\item Moreover, if $\Phi=(\ol{\Phi},\ul{\Phi})$ is an admissible map from  the bundle $\map{\tau}{E}{M}$ to the bundle $\map{\tau'}{E'}{M'}$, then we can identify the maps $\ol{\Phi}\equiv\prol[\ol{\Phi}]{\ul{\Phi}}$ and therefore $\Phi\equiv\prol[\Phi]{\Phi}$. 
\end{itemize}
\end{remark}

We finally mention that the composition of prolongation maps is the prolongations of the composition. Indeed, let $\Psi'$ be another bundle map from $\map{\pi'}{P'}{M'}$ to another bundle $\map{\pi''}{P''}{M''}$ and  $\Phi'$ be another admissible map from $\map{\tau'}{E'}{M'}$ to $\map{\tau''}{E''}{M''}$ both over the same base map. Since $\Phi$ and $\Phi'$ are admissible maps then so is $\Phi'\circ\Phi$, and thus we can define the prolongation of $\Psi'\circ\Psi$ with respect to $\Phi'\circ\Phi$. We have that  $\prol[\Phi'\circ\Phi]{(\Psi'\circ\Psi)}
=(\prol[\Phi']{\Psi'})\circ(\prol[\Phi]{\Psi})$.

\subsection*{Flow of a section} 
Every section of a Lie algebroid has an associated local flow composed of morphisms of Lie algebroids as we are going to explain next. In this context, by a flow on a bundle $\map{\tau}{E}{M}$ we mean a family of vector bundle maps $\Phi_s=(\ol{\Phi}_s,\ul{\Phi}_s)$ where $\ol{\Phi}_s$ is a flow (in the ordinary sense) on the manifold $E$ and $\ul{\Phi}_s$ is a flow (in the ordinary sense) on the manifold $M$.

\begin{theorem}
Let $\sigma$ be a section of a Lie algebroid $\map{\tau}{E}{M}$. There exists a local flow $\Phi_s$ on the bundle $\map{\tau}{E}{M}$ such that 
$$
d_\sigma\theta=\frac{d}{ds}\Phi_s\pb\theta\at{s=0},
$$
for any section $\theta$ of $\ext{E}$. For every fixed $s$, the map $\Phi_s$ is a morphism of Lie algebroids. The base flow $\map{\ul{\Phi}_s}{M}{M}$ is the flow (in the ordinary sense) of the vectorfield $\rho(\sigma)\in\vectorfields{M}$. 
\end{theorem}
\begin{proof}
Indeed, let $X_\sigma^C$ be the vector field on $E$ determined by its action on fiberwise linear functions,
$$
X_\sigma^C\hat{\theta}=\widehat{d_\sigma\theta},
$$
for every section $\theta$ of $E^*$, and which is known as the complete lift of $\sigma$ (see~\cite{GrUr,LMLA}). It follows that $X_\sigma^C\in\vectorfields{E}$ is projectable and projects to the vectorfield $\rho(\sigma)\in\vectorfields{M}$. Thus, if we consider the flow $\ol{\Phi}_s$ of $X_\sigma^C$ and the flow $\ul{\Phi}_s$ of $\rho(\sigma)$, we have that $\Phi_s=(\ol{\Phi}_s,\ul{\Phi}_s)$ is a flow on the bundle $\map{\tau}{E}{M}$. Moreover, since $X_\sigma^C$ is linear (maps linear functions into linear functions) we have that $\Phi_s$ are vector bundle maps. 

In order to prove the relation $d_\sigma\theta=\frac{d}{ds}\Phi_s\pb\theta\at{s=0}$, it is enough to prove it for 1-forms. If $\theta$ be a section of $E^*$, then for every $m\in M$ and $a\in E_m$ we have 
\begin{align*}
\pai{(d_\sigma\theta)_m}{a}
&=\widehat{d_\sigma\theta}(a)
=(X_\sigma^C\hat{\theta})(a)
=\frac{d}{ds}(\hat{\theta}\circ\ol{\Phi}_s)\at{s=0}(a)\\
&=\frac{d}{ds}\pai{\theta_{\ul{\Phi}_s(m)}}{\ol{\Phi}_s(a)}\at{s=0}
=\frac{d}{ds}\Phi_s\pb\theta\at{s=0}(a).
\end{align*}

Moreover, the maps $\Phi_s$ are morphisms of Lie algebroids. Indeed, from the relation $d\circ d_\sigma=d_\sigma\circ d$ we have that $\frac{d}{ds}(\Phi_s\pb d\theta-d\Phi_s\pb\theta)\at{s=0}=0$ for every section $\theta$ of $E^*$. It follows that $\Phi_s\pb d\theta=d\Phi_s\pb\theta$, so that $\Phi_s$ is a morphism.
\end{proof}

By duality, we have that 
$$
d_\sigma\zeta=\frac{d}{ds}\Phi_s\pb\zeta\at{s=0}
$$ 
for every section $\zeta$ of $E$.  Therefore we have a similar formula for the derivative of any tensorfield.

\begin{remark}
Alternatively one can describe the flow of a section in terms of Poisson geometry. A Lie algebroid structure is equivalent to a fiberwise linear Poisson structure on the dual $E^*$. The section $\sigma$ defines a linear function in $E^*$, and thus has an associated Hamiltonian vectorfield. The flow of this vectorfield is the transpose of the flow defined here. Since the flow of a Hamiltonian vectorfield leaves invariant the Poisson tensor its dual is a morphism of Lie algebroids. See~\cite{Fernandes} for the details.
\end{remark}

\begin{remark}
In the case of the standard Lie algebroid $E=TM$ we have that a section $\sigma$ is but a vector field on $M$. In this case, the flow defined by the section $\sigma$ is but $\Phi_s=(T\phi_s,\phi_s)$ where $\phi_s$ is the ordinary flow of the vectorfield $\sigma$.
\end{remark}

In the case of a projectable section of $\prol[E]{P}$ the flow can be expressed in terms of the flows of its components.

\begin{proposition}
\label{flow.projectable}
Let $Z$ be a projectable section of $\prol[E]{P}$ and let $\sigma\in\sec{E}$ the section to which it projects. Let $\Phi_s$ be the flow of the section $\sigma$. Let $\ol{\Psi}_s$ be the flow of the vectorfield $\rho^1(Z)\in\vectorfields{P}$, $\ul{\Psi}_s$ the flow of the vector field $\rho(\sigma)\in\vectorfields{M}$, and $\Psi_s=(\ol{\Psi}_s,\ul{\Psi}_s)$. The flow of $Z$ is $\prol[\Phi_s]{\Psi_s}$.
\end{proposition}
\begin{proof}
We first note that $\prol[\Phi_s]{\Psi_s}$ is well defined since the base maps of both $\Phi_s$ and $\Psi_s$ are equal to the flow of $\rho(\sigma)$. We have to prove that $d_Z\beta=\frac{d}{ds}(\prol[\Phi_s]{\Psi_s})\pb\beta\at{s=0}$, for every form $\beta\in \ext{\prol[E]{P}}$. But it is sufficient to prove it for 0-forms and 1-forms. For a function $F\in\cinfty{P}$ we have that 
$$
d_ZF=\rho^1(Z)F=\frac{d}{ds}F\circ\ol{\Psi}_s\at{s=0}=\frac{d}{ds}(\prol[\Phi_s]{\Psi_s})\pb F\at{s=0}.
$$
For 1-forms it is sufficient to prove the property for exact and basic 1-forms. For an exact form $\beta=dF$ we have
$$
d_ZdF=dd_ZF=d\frac{d}{ds}(\prol[\Phi_s]{\Psi_s})\pb F\at{s=0}
=\frac{d}{ds}d(\prol[\Phi_s]{\Psi_s})\pb F\at{s=0}
=\frac{d}{ds}(\prol[\Phi_s]{\Psi_s})\pb dF\at{s=0}.
$$
Finally for a basic form  $\beta=(\prol{\pi})\pb\theta$ we have
\begin{align*}
d_Z[(\prol{\pi})\pb\theta]
&=(\prol{\pi})\pb d_\sigma\theta
=(\prol{\pi})\pb\frac{d}{ds}\Phi_s\pb\theta\at{s=0}
=\frac{d}{ds}(\prol{\pi})\pb\Phi_s\pb\theta\at{s=0}\\
&=\frac{d}{ds}(\Phi_s\circ\prol{\pi})\pb\theta\at{s=0}
=\frac{d}{ds}(\prol{\pi}\circ\prol[\Phi_s]{\Psi_s})\pb\theta\at{s=0}
=\frac{d}{ds}(\prol[\Phi_s]{\Psi_s})\pb[(\prol{\pi})\pb\theta]\at{s=0},
\end{align*}
which ends the proof.
\end{proof}


\section{Jets}
\label{jetoids}

In this section we introduce our generalization of the space of jets, which will be the manifold where the Lagrangian of a field theory is defined. The key idea is to substitute the tangent bundle of a manifold by an appropriate Lie algebroid. To clarify the exposition we will introduce the concepts in an incremental manner, defining the concepts first on vector spaces, then on vector bundles and later introducing the anchor and finally the Lie algebroid structure. 

\subsection*{Splittings of exact sequences}
Let $F$ and $E$ be two finite dimensional real vector spaces and let  $\map{\pi}{E}{F}$ be a surjective linear map. If we denote by $K$ the kernel of $\pi$, we have the following exact sequence
$$
\seq 0-> K -{i}-> E -{\pi}-> F ->0,
$$
where $i$ is the canonical inclusion. Therefore we have the induced sequence
$$
\seq 0->\Lin(F,K) -{i\circ}-> \Lin(F,E) -{\pi\circ} -> \Lin(F,F)->0.
$$
For every linear map $\map{\zeta}{F}{F}$ we set
$$
\Lin^\zeta(F,E)=\set{\phi\in\Lin(F,E)}{\pi\circ\phi=\zeta}.
$$
It follows that $\Lin^\zeta(F,E)$ is an affine subspace of $\Lin(F,E)$ with underlying vector space $\Lin^0(F,E)=\Lin(F,K)$. For $\zeta=\id_F$ we have that $\Lin^\id(F,E)$ is the set of splittings of the original sequence and it is a model for all that affine spaces. The following diagram explain the situation
$$
\xymatrix{%
&0\ar[r]&\Lin(F,K)\ar[r]^{i\circ }&\Lin(F,E)\ar[r]^{\pi\circ}&\Lin(F,F)\ar[r]&0 \\
&&&\Lin^{\id}(F,E)\ar[u]^{}\ar[ur]_{constant=\id}&&
}.
$$

The vector space dual to $\Lin(F,E)$ is $\Lin(F,E)^*=\Lin(E,F)$, with the pairing
$$
\pai{p}{w}=\tr(p\circ w)=\tr(w\circ p),
$$
for $p\in\Lin(E,F)$ and $w\in\Lin(F,E)$. It follows that any element $p\in\Lin(E,F)$ can be identified with an element $\tilde{p}$ in the bidual space $\Lin(F,E)^*=\Lin(E,F)^{**}$, that is, with a linear function $\map{\tilde{p}}{\Lin(F,E)}{\Real}$ given by $\tilde{p}(w)=\tr(p\circ w)$.
By restricting to the affine space $\Lin^\id(F,E)$ we get affine functions
$\map{\hat{p}}{\Lin^\id(F,E)}{\Real}$ explicitly given by
$$
\hat{p}(\phi)=\tr(p\circ\phi).
$$
The linear function $\bar{p}$ associated to the affine function $\hat{p}$ is obtained by restriction of $\tilde{p}$ to $\Lin(F,K)$, that is, $\bar{p}(\psi)=\tr(p|_{K}\circ\Psi)$.

Obviously any linear map from $\Lin(F,E)$ to $\Real$ is of the form $\tilde{p}$ for a unique $p$. Also, every affine function is of the form $\hat{p}$, but notice that there can be many elements $p$ giving the same affine function. If $p, q\in\Lin(E,F)$ are such that $\hat{p}=\hat{q}$ then $p|_K=q|_K$ and $\tr(p-q)=0$, where $p-q$ is to be understood as an endomorphism of $F$.

Let us fix basis $\{\bar{e}_a\}$ in $F$ and $\{e_a,e_\alpha\}$ in $E$, adapted to $\pi$, that is, such that $\pi(e_a)=\bar{e}_a$ and $\pi(e_\alpha)=0$. Then we have the following expressions
\begin{align*}
&w=\bar{e}^b\otimes(y^a_be_a+y^\alpha_b e_\alpha) &&w\in\Lin(F,E)\\
&p=(p^a_be^b+p^a_\alpha e^\alpha)\otimes\bar{e}_a &&p\in\Lin(E,F)\\
&\phi=\bar{e}^a\otimes(e_a+y^\alpha_a e_\alpha)   &&\phi\in\Lin^\id(F,E)\\
&\psi=\bar{e}^a\otimes y^\alpha_a e_\alpha        &&\psi\in\Lin(F,K),
\end{align*}
and the corresponding linear and affine functions are
$$
\tilde{p}=p^a_by^b_a+p^a_\alpha y^\alpha_a
\qquand
\hat{p}=p^a_a+p^a_\alpha y^\alpha_a.
$$
The linear function associated to $\hat{p}$ is $\bar{p}=p^a_\alpha y^\alpha_a$. In particular we have that $\pi\in\Lin(E,F)$ and $\hat{\pi}=\tr(\id)=\dim F$ is constant. From this expression it is clear that two elements $p$ and $q$ define the same affine function if and only if $p^a_\alpha=q^a_\alpha$ and $p^a_b-q^a_b$ is trace free.

\subsection*{Vector bundles}
We consider now two vector bundles $\map{\prEM}{E}{M}$ and $\map{\prFN}{F}{N}$ and a surjective vector bundle map $\map{\prEF}{E}{F}$ over the map $\map{\prMN}{M}{N}$, as it is shown in the diagram
$$
\xymatrix{%
&E\ar[d]_{\prEM}\ar[r]^{\prEF}&F\ar[d]^{\prFN}\\
&M\ar[r]_{\prMN}&N
}
$$
Moreover, we will assume that $\prMN$ is a surjective submersion, so that  $\map{\prMN}{M}{N}$ is a smooth fiber bundle. By the symbol $\pi$ we mean the vector bundle map $\pi=(\prEF,\prMN)$. We will denote by $K\to M$ the kernel of the map $\pi$, which is a vector bundle over~$M$. Given a point $m\in M$, if we denote $n=\prMN(m)$, we have the following exact sequence
$$
\seq 0-> K_m -{}-> E_m -{}-> F_n ->0,
$$
and we can perform point by point similar constructions to those given above.

We will use the notations $\Lpi[m]\equiv\Lin(F_n,E_m)$, $\Jpi[m]\equiv\Lin^\id(F_n,E_m)$ and $\Vpi[m]\equiv\Lin(F_n,K_m)$, that is,
\begin{align*}
\Lpi[m]&=\set{\map{w}{F_n}{E_m}}{\text{$w$ is linear}}\\
\Jpi[m]&=\set{\phi\in\Lpi[m]}{\prEF\circ\phi=\id_{F_n}}\\
\Vpi[m]&=\set{\psi\in\Lpi[m]}{\prEF\circ\psi=0}.
\end{align*}
Therefore $\Lpi[m]$ is a vector space, $\Vpi[m]$ is a vector subspace of $\Lpi[m]$ and $\Jpi[m]$ is an affine subspace of $\Lpi[m]$ modeled on the vector space is $\Vpi[m]$. The dual vector space of $\Lpi[m]$ is $\dLpi[m]\equiv\Lin(E_m,F_n)$, that is,
$$
\dLpi[m]=\set{\map{p}{E_m}{F_n}}{\text{$p$ linear}}
$$

By taking the union, $\Lpi=\cup_{m\in M}\Lpi[m]$, $\Jpi=\cup_{m\in M}\Jpi[m]$ and $\Vpi=\cup_{m\in M}\Vpi[m]$, we get the vector bundle $\map{\prLM}{\Lpi}{M}$, the  vector subbundle $\map{\prVM}{\Vpi}{M}$ and the affine subbundle $\map{\prJM}{\Jpi}{M}$ modeled on the vector bundle $\Vpi$. The vector bundle dual to $\Lpi$ is $\dLpi=\cup_{m\in M}\dLpi[m]\to M$. We also consider the projection $\map{\ul{\pi}_1}{\Jpi}{N}$ given by composition $\prJN=\prMN\circ\prJM$.

Every section $\theta\in\Sec{\dLpi}$ defines the fiberwise linear function $\tilde{\theta}\in\cinfty{\Lpi}$ by $\tilde{\theta}(z)=\tr(\theta_m\circ z)$,
for $z\in\Lpi$ and where $m=\ul{\tilde{\pi}_{10}}(z)$. Its restriction to $\Jpi$ defines the affine function $\hat{\theta}\in\cinfty{\Jpi}$ by
$$
\hat{\theta}(\phi)=\tr(\theta_m\circ\phi),
$$
for $\phi\in\Jpi$ and where $m=\prJM(\phi)$.

\begin{remark}
In the standard case one has a bundle $\map{\nu}{M}{N}$ and then considers the tangent spaces at $m\in M$ and $n=\nu(m)\in N$ together with the differential of the projection $\map{T_m\nu}{T_mM}{T_nN}$. If we set $\pi=(T\nu,\nu)$, an element $\phi\in\Jpi[m]$ is a linear map $\map{\phi}{T_nN}{T_mM}$ such that $T_m\nu\circ\phi=\id_{T_nN}$. It follows that there exist sections $\varphi$ of $\nu$ such that $\varphi(n)=m$ and $T_n\varphi=\phi$. Thus, in the case of the tangent bundles $\map{\tau_N}{F\equiv TN}{N}$ and $\map{\tau_M}{E\equiv TM}{M}$ with  $\pi=(T\nu,\nu)$, an element $\phi\in\Jpi[m]$ can be identified with a  1-jet at the point $m$, i.e. an equivalence class of sections which has the same value and same first derivative at the point $m$. With the standard notations~\cite{Saunders}, we have that $J^1\nu\equiv\Jpi$. Obviously, this example leads our developments throughout this paper. 
\end{remark}

In view of this fact, an element of $\Jpi[m]$ will be simply called a \emph{jet} at the point $m\in M$ and accordingly the bundle $\Jpi$ is said to be the first \emph{jet bundle} of $\pi$.

\medskip

Local coordinates on $\Jpi$ are given as follows. We consider local coordinates $(\bar{x}^i)$ on $N$ and $(x^i,u^A)$ on $M$ adapted to the projection $\prMN$, that is $\bar{x}^i\circ\prMN=x^i$. We also consider local basis of sections $\{\bar{e}_a\}$ of $F$ and $\{e_a,e_\alpha\}$ of $E$ adapted to the projection $\prEF$, that is $\prEF\circ e_a=\bar{e}_a\circ\prMN$ and  $\prEF\circ e_\alpha=0$. In this way $\{e_\alpha\}$ is a base of sections of $K$. An element $w$ in $\Lpi[m]$ is of the form $w=(w_a^be_b+w_a^\alpha e_\alpha)\otimes \bar{e}^a$, and it is in $\Jpi[m]$ if and only if $w_a^b=\delta_a^b$, i.e.\ an element $\phi$ in $\Jpi$ is of the form $\phi=(e_a+\phi^\alpha_a e_\alpha)\otimes\bar{e}_a$. If we set $y^\alpha_a(\phi)=\phi^\alpha_a$, we have adapted local coordinates $(x^i,u^A,y^\alpha_a)$ on $\Jpi$. Similarly, an element $\psi\in\Vpi[m]$ is of the form $\psi=\psi^\alpha_a e_\alpha\otimes e^a$. If we set $y^\alpha_a(\psi)=\psi^\alpha_a$, we have adapted local coordinates $(x^i,u^A,y^\alpha_a)$ on $\Vpi$. As usual, we use the same name for the coordinates in an affine bundle and in the vector bundle on which it is modeled. A section of $\dLpi$ is of the form $\theta=(\theta^a_b(x)e^b+\theta^a_\alpha(x)e^\alpha)\otimes\bar{e}_a$, and the affine function defined by $\theta$ is
$$
\hat{\theta}=\theta^a_a(x)+\theta^a_\alpha(x)y^\alpha_a.
$$

\subsection*{Anchor}

We now assume that the given vector bundles are anchored vector bundles, that is, we have two vector bundle maps $\map{\rho_F}{F}{TN}$ and $\map{\rho_E}{E}{TM}$ over the identity in $N$ and $M$ respectively.

We will assume that the map $\pi$ is admissible, that is $d\pi\pb f=\pi\pb df$, where on the left $d$ stands for the exterior differential in $E$ while on the right $d$ stands for the exterior differential in $F$. Let us see what this condition means in terms of the components of the anchor.

\medskip

The anchor is determined by the differential of the coordinate functions $d\bar{x}^i=\bar{\rho}^i_a\bar{e}^a$ on $F$ and $dx^i=\rho^i_a e^a+\rho^i_\alpha e^\alpha$, $du^A=\rho^A_a e^a+\rho^A_\alpha e^\alpha$ on $E$.  Applying the admissibility condition to the coordinates functions $\bar{x}^i$ we get
$$
\bar{\rho}^i_a e^a
=\pi\pb(d\bar{x}^i)
=d(\pi\pb \bar{x}^i)
=dx^i
=\rho^i_a e^a+\rho^i_\alpha e^\alpha,
$$
where we have used that $\pi\pb\bar{e}^a=e^a$. Therefore  $\bar{\rho}^i_a=\rho^i_a$ (is a basic function) and $\rho^i_\alpha=0$. We have
$$
d\bar{x}^i=\rho^i_a\bar{e}^a
\qquand
\left\{\begin{array}{l}
dx^i=\rho^i_a e^a\\
du^A=\rho^A_a e^a+\rho^A_\alpha e^\alpha,
\end{array}\right.
$$
and equivalently
$$
\rho_F(\bar{e}_a)=\rho^i_a\pd{}{\bar{x}^i}
\qquand
\left\{\begin{array}{l}
\rho_E(e_a)=\rho^i_a\dpd{}{x^i}+\rho^A_a\dpd{}{u^A}\\[5pt]
\rho_E(e_\alpha)=\rho^A_\alpha\dpd{}{u^A},
\end{array}\right.
$$
with $\rho^i_a=\rho^i_a(x)$, $\rho^A_a=\rho^A_a(x,u)$ and $\rho^A_\alpha=\rho^A_\alpha(x,u)$.

\medskip

The anchor allows to define the concept of total derivative of a function with respect to a section. Given a section $\sigma\in\Sec{F}$, the total derivative of a function $f\in\cinfty{M}$ with respect to $\sigma$ is the function $\widehat{df\otimes\sigma}$, i.e. the affine function associated to $df\otimes\sigma\in\Sec{\dLpi}$. In particular, the total derivative with respect to the element $\bar{e}_a$ of the local basis of sections of $F$,  will be denoted by $\p{f}_{|a}$. In this way, if $\sigma=\sigma^a\bar{e}_a$ then $\widehat{df\otimes\sigma}=\p{f}_{|a}\sigma^a$,
where the coordinate expression of $\p{f}_{|a}$ is
$$
\p{f}_{|a}=\rho^i_a\pd{f}{x^i}+(\rho^A_a+\rho^A_\alpha y^\alpha_a)\pd{f}{u^A}.
$$
In particular, for the coordinate functions we have $\p{x}^i_{|a}=\rho^i_a$ and $\p{u}^A_{|a}=\rho^A_a+\rho^A_\alpha y^\alpha_a$, so that the chain rule $\p{f}_{|a}=\pd{f}{x^i}\p{x}^i_{|a}+\pd{f}{u^A}\p{u}^A_{|a}$ applies.

Notice that, for a basic function $f$ (i.e. the pullback of a function in the base $N$) we have that $\p{f}_{|a}=\rho^i_a\pd{f}{x^i}$ are just the (pullback of) components of $df$ in the basis $\{\bar{e}_a\}$.

\subsection*{Bracket}

Let us finally assume that we have Lie algebroid structures on $\map{\prFN}{F}{N}$ and on $\map{\prEM}{E}{M}$, and that the projection $\pi$ is a morphism of Lie algebroids. The kernel $K$ of $\pi$ is a $\pi$-ideal, that is, the bracket of two projectable sections is projectable, and the bracket of a projectable section with a section of $K$ is again a section of $K$.

This can be locally seen in the vanishing of some structure functions. The morphism condition $\pi\pb(d\theta)=d(\pi\pb\theta)$ applied to the 1-forms $\theta=\bar{e}^a$ gives
\begin{align*}
\pi\pb(d\bar{e}^a)
&=\pi\pb\left(-\frac{1}{2}\bar{C}^a_{bc}\bar{e}^b\wedge\bar{e}^c\right)
=-\frac{1}{2}\bar{C}^a_{bc}e^b\wedge e^c\\
d(\pi\pb\bar{e}^a)
&=-\frac{1}{2}C^a_{bc}e^b\wedge e^c
 -           C^a_{b\gamma}e^b\wedge e^\gamma
 -\frac{1}{2}C^a_{\beta\gamma}e^b\wedge e^\gamma,
\end{align*}
that is
$$
\bar{C}^a_{bc}=C^a_{bc}
\qquad
C^a_{b\gamma}=0
\qquand
C^a_{\beta\gamma}=0.
$$
It follows that $C^a_{bc}$ is a basic function and that the exterior differentials in $F$ and $E$ are determined by
$$
d\bar{e}^a=-\frac{1}{2}C^a_{bc}\bar{e}^b\wedge\bar{e}^c
\quand
\left\{\begin{array}{l}
de^a=-\frac{1}{2}C^a_{bc}e^b\wedge e^c\\[5pt]
de^\alpha=-\frac{1}{2}C^\alpha_{bc}e^b\wedge e^c
          -           C^\alpha_{b\gamma}e^b\wedge e^\gamma
          -\frac{1}{2}C^\alpha_{\beta\gamma}e^\beta\wedge e^\gamma.
 \end{array}\right.
$$
Equivalently, we have the following expressions for the various brackets
$$
[\bar{e}_a,\bar{e}_b]=C^c_{ab}\bar{e}_c
\qquand
\left\{\begin{aligned}
&[e_a,e_b]=C^\gamma_{ab}e_\gamma+C^c_{ab}e_c\\
&[e_a,e_\beta]=C^\gamma_{a\beta}e_\gamma \\
&[e_\alpha,e_\beta]=C^\gamma_{\alpha\beta}e_\gamma.
 \end{aligned}\right.
$$

\begin{remark}
It follows from our local calculations that one can avoid the extra bar in the notation as long as one is working with functions and sections of $F^*$ and $E^*$, but one has to be careful when working with sections of $F$ and $E$. More explicitly, one can omit the bar about coordinates $\bar{x}^i$ and forms $\bar{e}^a$ but we must use a bar when referring to the elements $\bar{e}_a$ of the basis of sections of $F$, in order to distinguish it from the sections $e_a$ of $E$. In what follows we will omit the bar over the coordinates $x^i$ but to avoid any confusion we will explicitly distinguish $\bar{e}^a$ and $\bar{e}_a$, from $e^a$ and $e_a$, respectively.
\end{remark}

\subsection*{Affine structure functions}
The structure functions can be combined to give some affine functions which contains part of the structure of Lie algebroid and that will naturally appear in our treatment of Lagrangian Field Theory.

We define the affine functions $Z^\alpha_{a\gamma}$ and $Z^\alpha_{ac}$ by
$$
Z^\alpha_{a\gamma}=\widehat{(d_{e_\gamma}e^\alpha)\otimes\bar{e}_a}
\qquand
Z^\alpha_{ac}=\widehat{(d_{e_c}e^\alpha)\otimes\bar{e}_a}.
$$
Explicitly, we have
$$
Z^\alpha_{a\gamma}=C^\alpha_{a\gamma}+C^\alpha_{\beta\gamma}y^\beta_a
\qquand
Z^\alpha_{ac}=C^\alpha_{ac}+C^\alpha_{\beta c}y^\beta_a.
$$
Indeed,
\begin{align*}
(d_{e_\gamma}e^\alpha)\otimes\bar{e}_a
&=i_{e_\gamma}\left(-\frac{1}{2}C^\alpha_{bc}e^b\wedge e^c-C^\alpha_{b\theta}e^b\wedge e^\theta-\frac{1}{2}C^\alpha_{\beta\theta}e^\beta\wedge e^\theta\right)
\otimes\bar{e}_a\\
&=(C^\alpha_{b\gamma}e^b+C^\alpha_{\beta\gamma}e^\beta)\otimes\bar{e}_a,
\end{align*}
and thus
$$
Z^\alpha_{a\gamma}
=\widehat{(C^\alpha_{b\gamma}e^b+C^\alpha_{\beta\gamma}e^\beta)
  \otimes\bar{e}_a}
=C^\alpha_{a\gamma}+C^\alpha_{\beta\gamma}y^\beta_a.
$$
Similarly
\begin{align*}
(d_{e_c}e^\alpha)\otimes\bar{e}_a
&=i_{e_c}\left(-\frac{1}{2}C^\alpha_{bd}e^b\wedge e^d-C^\alpha_{b\theta}e^b\wedge e^\theta-\frac{1}{2}C^\alpha_{\beta\theta}e^\beta\wedge e^\theta\right)
\otimes\bar{e}_a\\
&=(C^\alpha_{bc}e^b+C^\alpha_{\beta c}e^\beta)\otimes\bar{e}_a,
\end{align*}
and thus
$$
Z^\alpha_{ac}
=\widehat{(C^\alpha_{bc}e^b+C^\alpha_{\beta c}e^\beta)
  \otimes\bar{e}_a}
=C^\alpha_{ac}+C^\alpha_{\beta c}y^\beta_a.
$$
Notice that $\widehat{(d_{e_\gamma}e^b)\otimes\bar{e}_a}$=0 and that $\widehat{(d_{e_c}e^b)\otimes\bar{e}_a}=C^b_{ac}$. For completeness we will sometimes write $Z^b_{ac}=C^b_{ac}$ and $Z^b_{a\gamma}=0$.


\section{Repeated jets and second-order jets}
\label{2-jets}

The equations for a (standard) classical Field Theory are generally second-order partial differential equations. In geometric terms, this equations determine a submanifold of the manifold of second-order jets, also called 2-jets. Therefore, we need to define the analog of a second-order jet in this generalized setting. In order to do that we follow the identification of a 2-jet with a holonomic repeated first-order jet (see~\cite{Saunders}). Therefore, we need a Lie algebroid with base $\Jpi$ (which will be just prolongations of $\Jpi$ with respect to $E$) and a surjective morphism to $F$. Then we will consider jets for this morphism and we will select those who satisfy some holonomy condition. This holonomy condition is expressed in terms of the contact ideal, that will also be defined here.

\subsection*{Prolongation of $\Jpi$}

We consider the $E$-tangent bundle  $\prol[E]{\Jpi}$ to the jet manifold $\Jpi$. Recall that the fibre of this bundle at the point $\phi\in\Jpi[m]$ is
$$
\prol[E]{\Jpi}[\phi]=\set{(a,V)\in E_m\times T_\phi\Jpi}{\rho(a)=T\prJM(V)}.
$$ 
Local coordinates $(x^i,u^A)$ on $M$ and a local basis $\{e_a,e_\alpha\}$ of sections of $E$ determine local coordinates $(x^i,u^A,y^\alpha_a)$ on $\Jpi$ and a local basis $\{\X_a,\X_\alpha,\V_\alpha^a\}$ of sections of $\prol[E]{\Jpi}$ as was explained in section~\ref{preliminaries}. Explicitly, those sections are given by
\begin{align*}
&\X_a(\phi)=(\phi,e_a(m),\rho^i_a\pd{}{x^i}+\rho^A_a\pd{}{u^A})\\
&\X_\alpha(\phi)=(\phi,e_\alpha(m),\rho^A_\alpha\pd{}{u^A})\\
&\V^a_\alpha(\phi)=(\phi,0_m,\pd{}{y^\alpha_a}).
\end{align*}
With this settings, an element $Z=(\phi,a,V)\in\prol[E]{\Jpi}[\phi]$ 
can be written as $Z=a^b\X_b+a^\alpha\X_\alpha+V^\alpha_a\V^a_\alpha$, and the components $a$ and $V$ has the local expression $a=a^be_b+a^\alpha e_\alpha$ and 
$$
V=\rho^1(Z)
 =\rho^i_ba^b\pd{}{x^i}
  +(\rho^A_ba^b+\rho^A_\alpha a^\alpha)\pd{}{u^A}
  +V^\alpha_a\pd{}{y^\alpha_a}.
$$
The brackets of the elements in that basis are
\begin{align*}
&[\X_a,\X_b]=C^c_{ab}\X_c+C^\gamma_{ab}\X_\gamma
&&[\X_a,\X_\beta]=C^\gamma_{a\beta}\X_\gamma
&&[\X_\alpha,\X_\beta]=C^\gamma_{\alpha\beta}\X_\gamma\\
&[\X_a,\V^b_\beta]=0
&&[\X_\alpha,\V^b_\beta]=0
&&[\V^a_\alpha,\V^b_\beta]=0.
\end{align*}
The dual basis is to be denoted $\{\X^a,\X^\alpha,\V^\alpha_a\}$ and the exterior differential on $\prol[E]{\Jpi}$ is determined by the differential of the coordinate functions 
\begin{align*}
&dx^i=\rho^i_a \X^a\\
&du^A=\rho^A_a \X^a+\rho^A_\alpha \X^\alpha\\
&dy^\alpha_a=\V^\alpha_a,
\intertext{and the differential of the elements of the basis}
&d\X^a=-\frac{1}{2}C^a_{bc}\X^b\wedge\X^c\\
&d\X^\alpha=-\frac{1}{2}C^\alpha_{bc}\X^b\wedge\X^c
          -           C^\alpha_{b\gamma}\X^b\wedge\X^\gamma
          -\frac{1}{2}C^\alpha_{\beta\gamma}\X^\beta\wedge\X^\gamma\\
&d\V^\alpha_a=0.
\end{align*}

Finally, defining the map $\ol{\pi_{10}}=\prol{\ul{\pi_{10}}}$ we have that the bundle map $\pi_{10}=(\ol{\pi_{10}},\ul{\pi_{10}})$ is a surjective morphism from the Lie algebroid $\map{\tau^E_{\Jpi}}{\prol[E]{\Jpi}}{\Jpi}$ to the Lie algebroid $\map{\tau^E_M}{E}{M}$. We will also consider the projection $\prTJF(\phi,a,V)=\prEF(a)$, so that the bundle map $\pi_1=(\prTJF,\prJN)$ is a surjective morphism from the Lie algebroid $\map{\tau^E_{\Jpi}}{\prol[E]{\Jpi}}{\Jpi}$ to the Lie algebroid $\map{\tau^F_N}{F}{N}$. Notice that $\pi_1$ is just the composition  $\pi_1=\pi\circ\pi_{10}$.

\subsection*{Repeated jets}

Once we have defined the first jet manifold $\Jpi$, we can iterate the process by considering the Lie algebroids $\map{\prFN}{F}{N}$ and $\map{\tau^E_{\Jpi}}{\prol[E]{\Jpi}}{\Jpi}$, the projection $\pi_1$ 
$$
\xymatrix{%
&\prol[E]{\Jpi}\ar[d]_{\tau^E_{\Jpi}}\ar[r]^{\prTJF}&F\ar[d]^{\prFN}\\
&\Jpi\ar[r]_{\prJN}&N
}
$$
and the set of jets of the above diagram which is the manifold $\JJpi$ fibred over ${\Jpi}$. An element $\psi$ of $\JJpi$ at the point $\phi\in\Jpi$ is thus a linear map  $\map{\psi}{F_n}{\CMcal{T}^{E}_\phi\Jpi}$ such that $\prTJF\circ\psi=\id_{F_n}$, where $n=\prJN(\phi)$. Let us see the explicit form of $\psi$. 

\begin{proposition}
An element $\psi\in\JJpi$ is of the form $(\phi,\zeta,V)$ for $\phi$, $\zeta\in\Jpi$ with $\prJM(\phi)=\prJM(\zeta)$ and $\map{V}{F_{\prJN(\phi)}}{T_\phi\Jpi}$ a linear map satisfying $T\prJM\circ V=\rho_E\circ\zeta$. 
\end{proposition}
\begin{proof}
Indeed, $\psi$ is a linear map $\map{\psi}{F_n}{\CMcal{T}^{E}_\phi\Jpi}$, so that it is of the form $\psi(b)=(\phi,\zeta(b),V(b))$ for some linear maps $\map{\zeta}{F_n}{E_m}$ and $\map{V}{F_n}{T_\phi\Jpi}$. The condition $\prTJF\circ\psi=\id_{F_n}$ is just $\prEF(\zeta(b))=b$, for every $b\in F_n$, i.e. $\prEF\circ\zeta=\id_{F_n}$, which is just to say that $\zeta$ is an element of $\Jpi$. Moreover, $\psi(b)$ is an element of $\prol[E]{\Jpi}$, so that we have $T\prJM(V(b))=\rho_E(\zeta(b))$, for every $b\in F_n$. In other words $T\prJM\circ V=\rho_E\circ\zeta$.
\end{proof}

Local coordinates $(x^i,u^A)$  and a local base of sections $\{e_a,e_\alpha\}$ as before, provide natural local coordinates for $\JJpi$. The three components of an element $\psi=(\phi,\zeta,V)\in\JJpi$ are locally of the form
$$
\phi=(e_a+y^\alpha_ae_\alpha)\otimes\bar{e}^a,
\qquad\qquad
\zeta=(e_a+z^\alpha_ae_\alpha)\otimes\bar{e}^a,
$$
and
$$
V=\left(
 \rho^i_a\pd{}{x^i}
+(\rho^A_a+\rho^A_\alpha z^\alpha_a)\pd{}{u^A}
+y^\beta_{ba}\pd{}{y^\beta_b}
\right)\otimes\bar{e}^a.
$$
Therefore, $(x^i,u^A,y^\alpha_a,z^\alpha_a,y^\alpha_{ba})$ are local coordinates for $\psi$. In terms of this coordinates and the associated local basis of $\prol[E]{\Jpi}$ we have the local expression
$$
\psi=(\X_a+z^\alpha_a\X_\alpha+y^\alpha_{ba}\V^b_\alpha)\otimes\bar{e}_a.
$$

\subsection*{Contact forms}
Contact forms are 1-forms on $\prol[E]{\Jpi}$ which satisfies a property of verticality, as it is explained in what follows. A jet $\phi\in\Jpi[m]$, being a splitting of an exact sequence, determines two complementary projectors $\map{h_{\phi}}{E_m}{E_m}$ and $\map{v_{\phi}}{E_m}{E_m}$ given by
$$
h_\phi(a)=\phi(\prEF(a))
\qquand
v_\phi(a)=a-\phi(\prEF(a)).
$$
An element $Z=(\phi,a, V)\in\prol[E]{\Jpi}$ is said to be horizontal if $a$ is horizontal with respect to $\phi$, that is, $v_\phi(a)=0$. In coordinates a horizontal element is of the form 
$$
Z=a^b(\X_b|_\phi+y^\beta_b\X_\beta|_\phi)+V^\beta_b\V^b_\beta|_\phi,
$$
where $y^\alpha_a$ are the coordinates of~$\phi$.
Let $(\prol[E]{\Jpi})^*$ the vector bundle dual to $\prol[E]{\Jpi}$. An element $\mu\in(\prol[E]{\Jpi}[\phi])^*$ is said to be vertical if it vanishes on every horizontal element at $\phi$. It follows that  a vertical element is of the form $\mu=\mu_\alpha(\X^\alpha|_\phi-y^\alpha_a\X^a|_\phi)$, where $y^\alpha_a$ are the coordinates of~$\phi$. 

A \emph{contact 1-form} is a section of $(\prol[E]{\Jpi})^*$ which is vertical at every point. It follows that it is a semibasic form and that the set of contact 1-forms is spanned by the forms $\theta^\alpha=\X^\alpha-y^\alpha_a\X^a$.

Every contact 1-form can be obtained as follows. Given a section $\alpha$ of $\prJM^*(E^*)$ (a section of $E^*$ along $\prJM$) we define the contact 1-form $\breve{\alpha}$ by means of $\pai{\breve{\alpha}}{(\phi,a,V)}=\pai{\alpha_\phi}{v_\phi(a)}$, for $(\phi,a,V)\in\Jpi$. 

The $\ext{(\prol[E]{\Jpi})}$-module generated by contact 1-forms is said to be the \emph{contact module} and will be denoted by $\mathcal{M}^c$. Since $\prol[E]{\Jpi}$ is a Lie algebroid, we have a differential on it, and we can consider the differential ideal generated by contact 1-forms. This ideal is said to be the \emph{contact ideal} and will be denoted $\CMcal{I}^c$. The contact ideal and the contact module are different (as sets), that is the contact ideal \textsc{is not} generated by contact 1-forms. Indeed, an easy calculation shows that
\begin{gather*}
d\theta^\alpha
+\frac{1}{2}C^\alpha_{\beta\gamma}\theta^\beta\wedge\theta^\gamma
+Z^\alpha_{b\gamma}\X^b\wedge\theta^\gamma
=
{}\hbox to 6cm{\hss}\\
\hbox to 2cm{\hss}{}=\X^a\wedge\V^\alpha_a+\frac{1}{2}\left(
y^\alpha_a C^a_{bc}
 - C^\alpha_{\beta\gamma}y^\beta_by^\gamma_c
 - C^\alpha_{b\gamma}y^\gamma_c+C^\alpha_{c\gamma}y^\gamma_b
 + C^\alpha_{bc}
\right)\X^b\wedge\X^c,
\end{gather*}
which cannot be written in the form $A^\alpha_\beta\wedge\theta^\beta$.

\subsection*{Second-order jets}

We define here the analog of the second order jet bundle as a submanifold of the repeated first jet bundle, in the same way as in the standard theory the manifold $J^2\nu$ can be considered as the submanifold of $J^1(J^1\nu)$ of holonomic jets.

\begin{definition}
A jet $\psi\in\JJpi[\phi]$ is said to be \emph{semiholonomic} if $\psi\pb\theta=0$ for every element $\theta$ in the contact module $\mathcal{M}^c$ at $\phi$. The set of all semiholonomic jets will be denoted $\widehat{\JJpi}$, that is,
$$
\widehat{\JJpi}=\set{\psi\in\JJpi}{\psi\pb\theta=0\ \text{for every}\ \theta\in\CMcal{M}^c}.
$$
\end{definition}

\begin{proposition}
An element $\psi\in\JJpi[\phi]$ is semiholonomic if and only if it is of the form  $\psi=(\phi,\phi,V)$. In other words if and only if it satisfies $\prol{\prJM}\circ\psi=\tau^E_{\Jpi}(\psi)$.
\end{proposition}
\begin{proof}
Indeed, this is equivalent to $v_\phi\circ\zeta=0$ which implies $\zeta=\phi$. A proof in coordinates is as follows. If the coordinates of $\phi$ are $(x^i,u^A,y^\alpha_a)$, and the local expression of $\psi\in\JJpi[\phi]$ is $\psi=(\X_a+z^\alpha_a\X_\alpha+y^\alpha_{ba}\V^b_\alpha)\otimes\bar{e}_a$, then 
$$
\psi\pb\theta^\alpha
=\psi\pb(\X^\alpha-y^\alpha_a\X^a)
=(z^\alpha_a-y^\alpha_a)\bar{e}^\alpha,
$$
so that $\psi\pb\theta^\alpha=0$ if and only if $z^\alpha_a=y^\alpha_a$ which is equivalent to $\phi=\zeta$.
\end{proof}

In the local coordinate system $(x^i,u^A,y^\alpha_a,z^\alpha_a,y^\alpha_{ab})$ on $\JJpi$, the coordinates of a semiholonomic jet are of the form $(x^i,u^A,y^\alpha_a,y^\alpha_a,y^\alpha_{ab})$, and thus the set $\widehat{\JJpi}$ of semiholonomic jets is a submanifold  of $\JJpi$ in which we have local coordinates $(x^i,u^A,y^\alpha_a,y^\alpha_{ab})$.

\begin{definition}
A jet $\psi\in\JJpi[\phi]$ is said to be \emph{holonomic} if $\psi\pb\theta=0$ for every element $\theta$ in the contact ideal $\mathcal{I}^c$ at $\phi$. The set of all holonomic jets will be denoted $\Jtpi$, that is,
$$
\Jtpi=\set{\psi\in\JJpi}{\psi\pb\theta=0\ \text{for every}\ \theta\in\CMcal{I}^c}.
$$
\end{definition}

Holonomic jets play the role of second-order jets in the case of the standard theory of jet bundles, and we frequently refer to a holonomic jet as a second order jet or simply as a 2-jet.

From the definition it follows that a holonomic jet is necessarily semiholonomic, but notice that not every semiholonomic jet is holonomic. 

\begin{proposition}
A semiholonomic jet $\psi$ is holonomic if its coordinates satisfy the equations $\mathcal{M}^\alpha_{ab}=0$, where the local functions $\mathcal{M}^\alpha_{ab}$ are defined by 
$$
\mathcal{M}^\gamma_{ab}\equiv
y^\gamma_{ab}-y^\gamma_{ba}
+C^\gamma_{b\alpha}y^\alpha_a-C^\gamma_{a\beta}y^\beta_b
-C^\gamma_{\alpha\beta}y^\alpha_ay^\beta_b+y^\gamma_c C^c_{ab}
+C^\gamma_{ab}.
$$
\end{proposition}
\begin{proof}
Indeed, from the above expression of $d\theta^\alpha$ we have
\begin{align*}
\psi\pb d\theta^\alpha
&=\psi\pb(\X^b\wedge\V^\alpha_b+\frac{1}{2}\left(
y^\alpha_a C^a_{bc}
 - C^\alpha_{\beta\gamma}y^\beta_by^\gamma_c
 - C^\alpha_{b\gamma}y^\gamma_c+C^\alpha_{c\gamma}y^\gamma_b
 + C^\alpha_{bc}
\right)\X^b\wedge\X^c)\\
&=\frac{1}{2}\left(
y^\alpha_{bc}-y^\alpha_{cb}+
y^\alpha_a C^a_{bc}
 - C^\alpha_{\beta\gamma}y^\beta_by^\gamma_c
 - C^\alpha_{b\gamma}y^\gamma_c+C^\alpha_{c\gamma}y^\gamma_b
 + C^\alpha_{bc}
\right)\bar{e}^b\wedge\bar{e}^c\\
&=\frac{1}{2}\mathcal{M}^\alpha_{bc}\bar{e}^b\wedge\bar{e}^c
\end{align*}
from where the result immediately follows.
\end{proof}

\begin{remark}
The above condition is to be interpreted point by point. That is, we take the value of any contact form at the point $\phi$ and then take the pullback. We did not write explicitly the point $\phi$ in order to simplify the reading of the formulas.

The condition $\mathcal{M}^\alpha_{bc}=0$ admits a geometrical interpretation as the condition for $\psi$ to be the jet of a morphism of Lie algebroids (see section~\ref{morphisms}). Therefore we refer to it as the \emph{infinitesimal morphism} condition. 
\end{remark}

\begin{remark}
In the standard case where $E=TM$ and $F=TN$ with coordinate basis of vectorfields the above conditions reduce to $u^A_{ij}=u^A_{ji}$, which is the pointwise condition for the matrix $u^A_{ij}$ to be the Hessian of some functions $u^A$.
\end{remark}

From the above proposition we see that $\Jtpi$ is a submanifold of $\JJpi$ and that local coordinates for $\Jtpi$ are $(x^i,u^A,y^\alpha_a,y^\alpha_{ab})$ for $a\leq b$. The coordinates $y^\alpha_{ab}$ for $a>b$ are determined by the equations $\mathcal{M}^\alpha_{ab}=0$.


\section{Morphisms and admissible maps}
\label{morphisms}

By a section of $\pi$ we mean a vector bundle map $\Phi$ such that $\pi\circ\Phi=\id$, which explicitly means $\prMN\circ\ul{\Phi}=\id_N$ and $\prEF\circ\ol{\Phi}=\id_F$ (In other words we consider only linear sections of $\prEF$). The set of sections of $\pi$ will be denoted by $\sec{\pi}$. The set of those sections of $\pi$ which are admissible as a map between anchored vector bundles will be denoted by $\Adm{\pi}$, and the set of those sections of $\pi$ which are a morphism of Lie algebroids will be denoted by $\Mor{\pi}$. Clearly $\Sec{\pi}\subset\Adm{\pi}\subset\Mor{\pi}$. We will find in this section local conditions for a section to be an admissible map between anchored vector bundles and local conditions for a section to be a morphism of Lie algebroids.

Taking adapted local coordinates $(x^i,u^A)$ on $M$, the map $\ul{\Phi}$ has the expression $\ul{\Phi}(x^i)=(x^i,\phi^A(x))$. If we moreover take an adapted basis $\{e_a,e^\alpha\}$ of local sections of $E$, then the expression of $\ol{\Phi}$ is given by $\ol{\Phi}(\bar{e}_a)=e^a+\phi^\alpha_ae_\alpha$, so that the map $\Phi$ is determined by the functions $\bigl(\phi^A(x),\phi^\alpha_a(x)\bigr)$ locally defined on $N$. The pullback of the coordinate functions is $\Phi\pb x^i=x^i$ and $\Phi\pb u^A=\phi^A$ and the pullback of dual basis is $\Phi\pb e^a=\bar{e}^a$, and $\Phi\pb e^\alpha=\phi^\alpha_a \bar{e}^a$

Let us see what the admissibility condition $\Phi\pb(df)=d(\Phi\pb f)$ means for this maps. As before, we impose this condition to the coordinate functions. Taking $f=x^i$ we get an identity (i.e. no new condition arises), and taking $f=u^A$ we get
$$
d\phi^A
=d(\Phi\pb u^A)
=\Phi\pb(du^A)
=\Phi\pb(\rho^A_ae^a+\rho^A_\alpha e^\alpha)
=[(\rho^A_a\circ\ul{\Phi})+(\rho^A_\alpha\circ\ul{\Phi})\phi^\alpha_a] \, \bar{e}^\alpha,
$$
from where we get that $\Phi$ is an admissible map if and only if
$$
\rho^i_a\pd{\phi^A}{x^i}
=(\rho^A_a\circ\ul{\Phi})+(\rho^A_\alpha\circ\ul{\Phi})\phi^\alpha_a.
$$
In order to simplify the writing, we can omit the composition with $\ul{\Phi}$, since we know that this is an equation to be satisfied at the point $m=\ul{\Phi}(n)=(x^i,\phi^A(x))$ for every $n\in N$. (This is but the usual practice). With this convention, the above equation is written as
$$
\rho^i_a\pd{u^A}{x^i}=\rho^A_a+\rho^A_\alpha y^\alpha_a.
$$

Let us now see what the condition of being a morphism means in coordinates. If we impose $\Phi\pb de^a=d\Phi\pb e^a$ we get an identity, so that we just have to impose $\Phi\pb de^\alpha=d\Phi\pb e^\alpha$. On one hand we have
$$
d(\Phi\pb e^\alpha)
=d(\phi^\alpha_a\bar{e}^a)
=\frac{1}{2}\left(
  \rho^i_b\pd{\phi^\alpha_c}{x^i}-\rho^i_c\pd{\phi^\alpha_b}{x^i}
 -\phi^\alpha_a C^a_{bc}
  \right)\bar{e}^b\wedge\bar{e}^c
$$
and on the other hand
\begin{align*}
\Phi\pb d(e^\alpha)
&=-\Phi\pb\left(
   \frac{1}{2}C^\alpha_{\beta\gamma}e^\beta\wedge e^\gamma
  -C^\alpha_{b\gamma}e^b\wedge e^\gamma
  -\frac{1}{2}C^\alpha_{bc}e^b\wedge e^c
     \right)\\
&=-\frac{1}{2}\left(
   C^\alpha_{\beta\gamma}\phi^\beta_b\phi^\gamma_c
  +C^\alpha_{b\gamma}\phi^\gamma_c-C^\alpha_{c\gamma}\phi^\gamma_b
  -C^\alpha_{bc}\right)\bar{e}^b\wedge\bar{e}^c
\end{align*}
Thus, the bundle map $\Phi$ is a morphism if and only if it satisfies
$$
\rho^i_b\pd{\phi^\alpha_c}{x^i}-\rho^i_c\pd{\phi^\alpha_b}{x^i}
 -\phi^\alpha_a C^a_{bc}
+
   C^\alpha_{\beta\gamma}\phi^\beta_b\phi^\gamma_c
  +C^\alpha_{b\gamma}\phi^\gamma_c-C^\alpha_{c\gamma}\phi^\gamma_b
=C^\alpha_{bc},
$$
in addition to the admissibility condition. As before, in this equation is to be satisfied at every point $m=\ul{\Phi}(n)$ in the image of $\ul\Phi$.

Taking into account our notation $\p{f}_{|a}=\rho^i_a\pd{f}{x^i}$ for a function $f\in\cinfty{N}$ we can write the above expressions as
\begin{gather*}
\p{u}^A_{|a}=\rho^A_a+\rho^A_\alpha y^\alpha_a\\
 \p{y}^\alpha_{c|b}-\p{y}^\alpha_{b|c}
 +C^\alpha_{b\gamma}y^\gamma_c-C^\alpha_{c\gamma}y^\gamma_b
 +C^\alpha_{\beta\gamma}y^\beta_by^\gamma_c-y^\alpha_a C^a_{bc}
=C^\alpha_{bc}.
\end{gather*}
where we recall that this equations are to be satisfied at every point $m=\ul{\Phi}(n)$. Notice the similarity between the above expression and the infinitesimal morphism condition $\mathcal{M}^\alpha_{bc}=0$. The following subsections will make clear this relation.

\begin{remark}
In the standard case where $E=TM$ and $F=TN$ with coordinate basis of vectorfields the above morphism conditions reduce to
$$
y^A_i=\pd{u^A}{x^i}
\qquand
\pd{y^A_i}{x^j}=\pd{y^A_j}{x^i}.
$$
It follows that every admissible map is a morphism, and moreover is the tangent map of a section of $\nu$. Therefore, in the standard case, by considering morphisms we are just considering 1-jet prolongation of sections of the base bundle.
\end{remark}

\subsection*{Sections of $\ul{\pi_1}$}
At this point we want to make clear the equivalence between sections of $\map{\prJN}{\Jpi}{N}$ and vector bundle maps which are sections of  $\pi$. 

Indeed, let $\Phi\in\sec{\pi}$ be a section of $\pi$. For every point $n\in N$, the restriction of $\ol{\Phi}$ to the fibre $F_n$ is a map from $F_n$ to the fibre $E_m$, where $m=\ul{\Phi}(n)$. This map $\map{\ol{\Phi}_n\equiv\ol{\Phi}\big|_{F_n}}{F_n}{E_m}$ satisfies $\prEF\circ\ol{\Phi}_n=\id_{F_n}$ and thus $\ol{\Phi}_n\in\Jpi[m]$. In this way, we have defined a map $\map{\check{\ol{\Phi}}}{N}{\Jpi}$ given by $\check{\ol{\Phi}}(n)=\ol{\Phi}_n$. The map $\check{\ol{\Phi}}$ is a section of $\prJN$; indeed, for every $n\in N$,
$$
(\prJN\circ\check{\ol{\Phi}})(n)
=\prJN(\check{\ol{\Phi}}(n))
=\prMN\bigl(\prJM(\ol{\Phi}_n)\bigr)=\prMN(m)
=n.
$$ 

Conversely, given a section $\map{\Psi}{N}{\Jpi}$ of $\prJN$, we define the map $\map{\ul{\Phi}}{N}{M}$ by $\ul{\Phi}=\prJM\circ\Psi$. Thus $\ul{\Phi}$ is a section of $\prMN$ and for every $n\in N$ we have that $\Psi(n)\in\Jpi[\ul\Phi(n)]$. Therefore $\Psi(n)$, being a jet, is a map $\map{\Psi(n)}{F_n}{E_{\ul\Phi(n)}}$ such that $\prEF\circ\ul{\Phi}=\id_{F_n}$. Consider now the map $\map{\ol{\Phi}}{F}{E}$ given by $\ol{\Phi}(b)=\Psi(n)(b)$ for every  $b\in F$, where $n=\prFN(b)$.  In this way we have that $\Phi=(\ol{\Phi},\ul{\Phi})$ is a linear bundle map and  
$$
(\prEF\circ\ol{\Phi})(b)=\prEF(\Psi(n)(b))=b
$$
for every $b\in F$ and where $n=\prFN(b)$, so that $\Phi$ is a section of $\pi$. By construction, it is clear that $\ol{\Phi}\big|_{F_n}=\Psi(n)$ for every $n\in N$, so that $\check{\ol{\Phi}}=\Psi$. Thus we have proved the following result.

\begin{proposition}
There is a one to one correspondence between the set $\sec{\pi}$ of vector bundle maps which are sections of $\pi$ and the set of sections of $\ul{\pi_1}$.
\end{proposition}

Moreover, $\check{\ol{\Phi}}$ is a fibred map over $\ul{\Phi}$, i.e. $\ul{\pi_{10}}\circ\check{\ol{\Phi}}=\ul{\Phi}$. We will denote by $\check{\Phi}$ the bundle map $\check{\Phi}=(\check{\ol{\Phi}},\ul{\Phi})$ (from the bundle $\map{\id}{N}{N}$ to the bundle $\map{\ul{\pi_{10}}}{\Jpi}{M}$).

The above equivalence is important since we will look for equations to be satisfied by a vector bundle map $\Phi$, and this will be reformulated as an equation on $\Jpi$, which is to be satisfied by the associated map $\check{\Phi}$. Moreover, we will impose to $\Phi$ to be a morphism of Lie algebroids, so that we need a condition expressed in terms of jets equivalent to the admissibility and the morphism condition.

\subsection*{Jet prolongation of sections of $\pi$}

Let $\Phi=(\ol{\Phi},\ul{\Phi})$ be a section of $\pi$ which is admissible as a map between anchored vector bundles, i.e. $\Phi\in\Adm{\pi}$, and we consider the associated map $\map{\check{\Phi}}{N}{\Jpi}$ of $\prJN$. Since the base map of both $\Phi$ and $\check{\Phi}$ is $\ul{\Phi}$ and $\Phi$ is admissible, we can construct the prolongation map $\map{\prol[\ol{\Phi}]{\check{\ol{\Phi}}}}{F\equiv\prol[F]{N}}{\prol[E]{\Jpi}}$, and hence the bundle map $\Phi\spi=\prol[\Phi]{\check{\Phi}}=(\prol[\ol{\Phi}]{\check{\Phi}},\check{\ol{\Phi}})$ from $\map{\tau^F_N}{F}{N}$ to $\map{\tau^E_{\Jpi}}{\prol[E]{\Jpi}}{\Jpi}$. This map is a section of $\pi_1$ and moreover it is admissible $\Phi\spi\in\Adm{\pi_1}$. The section of $\map{\ul{(\pi_1)_1}}{\JJpi}{N}$ associated to $\Phi\spi$ will be denoted $\check{\Phi}\spi$ . 

\begin{definition}
The map $\Phi\spi$ is said to be the first  \emph{jet prolongation} of $\Phi$, and the section $\check{\Phi}\spi$ is said to be the first jet prolongation of the section $\check{\Phi}$.
\end{definition}

\begin{proposition}
For every admissible section $\Phi\in\Adm{\pi}$ the section $\check{\Phi}\spi$ is semiholonomic.
\end{proposition}
\begin{proof}
Indeed, the explicit expression of $\Phi\spi$ is
$$
\Phi\spi_n(b)=\prol[\ol{\Phi}]{\check{\ol{\Phi}}}(n,b,\rho_F(b))
=\bigl(\check{\Phi}(n),\ol{\Phi}(b),T\check{\Phi}(\rho_F(b))\bigr)
=\bigl(\Phi_n,\Phi_n(b),T\check{\Phi}(\rho_F(b))\bigr).
$$
so that $\check{\Phi}\spi(n)=(\Phi_n,\Phi_n,T\check{\Phi}\circ\rho_F\bigr)$ is semiholonomic for every $n\in N$.
\end{proof}

In local coordinates $(x^i,u^A,y^\alpha_a,z^\alpha_a,y^\alpha_{ab})$ in $\JJpi$, if the expression of $\check{\Phi}$ is  $\check{\Phi}(x)=(x^i,\phi^A,\phi^\alpha_a)$, then the expression of $\check{\Phi}\spi$ is $\check{\Phi}\spi(x)= (x^i,\phi^A,\phi^\alpha_a,\phi^\alpha_a,\p{\phi}^\beta_{b|a})$. Indeed, from the general expression of the prolongation of a map with respect to a morphism we have 
\begin{align*}
&\Phi\spi{}\pb\X^a
 =(\prol[\Phi]{\check{\Phi}})\pb\X^a
 =\bar{e}^a\\
&\Phi\spi{}\pb\X^\alpha
 =(\prol[\Phi]{\check{\Phi}})\pb\X^\alpha
 =\phi^\alpha_a\bar{e}^a\\
&\Phi\spi{}\pb\V^\beta_b
 =(\prol[\Phi]{\check{\Phi}})\pb\V^\beta_b
 =d\phi^\beta_b=\p{\phi}^\beta_{b|a}\bar{e}^a,
\end{align*} 
so that 
$$
\Phi\spi=(\X_a+\phi^\alpha_a\X_\alpha+\p\phi^\beta_{b|a}\V^\beta_b)
\otimes\bar{e}^a,
$$
which is equivalent to the above coordinate expression.

From this local result we see that, if $\Phi$ is an admissible section of $\pi$, then $\Phi$ is a morphism if and only if $\check{\Phi}\spi$ is holonomic, i.e. takes values in $\Jtpi$. In the following subsection we will perform a more detailed study of this fact.

\subsection*{Morphisms and holonomic jets}

The equations for a classical Field Theory as they arise from variational calculus are equations for a morphism of Lie algebroids~\cite{DGA04}. On the other hand, the equations that we are going to set are equations for jets. Therefore, it is necessary to find conditions that select jets of  morphism among the set of jet of general maps. In this section we will find such condition in terms of holonomy.

\begin{proposition}\label{admissible-prolongation}
Let $\Psi\in\sec{\pi_1}$ be such that the associated map $\check{\Psi}$ is a semiholonomic section and let $\check{\Phi}$ be the section of $\ul{\pi_1}$ to which it projects. Then 
\begin{enumerate}
\item The bundle map $\Psi$ is admissible  if and only if $\Phi$ is admissible and $\Psi=\Phi\spi$. 
\item The bundle map $\Psi$ is a morphism of Lie algebroids if and only if $\Psi=\Phi\spi$ and $\Phi$ is a morphism of Lie algebroids.
\end{enumerate}
\end{proposition}
\begin{proof}
In local adapted coordinates $(x^i,u^A,y^\alpha_a,z^\alpha_a,y^\alpha_{ab})$ on $\JJpi$ we have that the section $\check{\Psi}$ associated to $\Psi$ is 
$\check{\Psi}(x)=(x^i,\phi^A(x),\phi^\alpha_a(x),\psi^\beta_a,\psi^\beta_{ba})$
so that the action of $\Psi$ by pullback is 
\begin{align*}
&\Psi\pb x^i=x^i
&&\Psi\pb\X^a=\bar {e}^a\\
&\Psi\pb u^A=\phi^A
&&\Psi\pb\X^\beta=\psi^\beta_a\bar {e}^a\\
&\Psi\pb y^\alpha_a=\phi^\alpha_a
&&\Psi\pb\V^\beta_b=\psi^\beta_{ba}\bar {e}^a.
\end{align*}

Let us see the conditions for $\Psi$ to be an admissible map between anchored bundles. We impose $\Psi\pb df=d\Psi\pb f$  for the coordinate functions. For $f=x^i$ we get an identity. For $f=u^A$ we have
$$
\Psi\pb du^A-d\Psi\pb u^A
=\Psi\pb(\rho^A_a\X^a+\rho^A_\alpha\X^\alpha)-d\phi^A
=\left(\rho^A_a+\rho^A_\alpha\psi^\alpha_a-\rho^i_a\pd{\phi^A}{x^i}\right)\bar{e}^a.
$$
Finally for $f=y^\beta_b$ we have
$$
\Psi\pb dy^\beta_b-d\Psi\pb y^\beta_b
=\Psi\pb\V^\beta_b-d\phi^\beta_b
=\left(\psi^\beta_{ba}-\rho^i_a\pd{\phi^\beta_b}{x^i}\right)\bar{e}^a
$$
Therefore, the map $\Psi$ is admissible if and only if
\begin{enumerate}
\item[(a)] $\p{\phi}^A_{|a}=\rho^A_a+\rho^A_\alpha\psi^\alpha_a$, and
\item[(b)] $\p{\phi}^\beta_{b|a}=\psi^\beta_{ba}$.
\end{enumerate}

We now impose $\Psi\pb d\theta=d\Psi\pb\theta$ for $\theta$ an element in the dual basis $\{\X^a,\X^\beta,\V^\beta_b\}$. For $\theta=\X^a$ we get an identity. For $\theta=\X^\alpha$ we have
\begin{align*}
\Psi\pb d\X^\alpha
&=\Psi\pb\left(
 -\frac{1}{2}C^\alpha_{\beta\gamma}\X^\beta\wedge\X^\gamma
 -C^\alpha_{b\gamma}\X^b\wedge\X^\gamma
 -\frac{1}{2}C^\alpha_{bc}\X^b\wedge\X^c\right)\\
&=-\frac{1}{2}\left(
  C^\alpha_{\beta\gamma}\psi^\beta_b\psi^\gamma_c
 +C^\alpha_{b\gamma}\psi^\gamma_c
 -C^\alpha_{c\beta}\psi^\beta_b
 -C^\alpha_{bc}\right)
 \bar{e}^b\wedge\bar{e}^c,
\intertext{and} 
d\Psi\pb \X^\alpha
&=d(\psi^\alpha_a\bar{e}^a)
=-\frac{1}{2}\left(
  \rho^i_c\pd{\psi^\alpha_b}{x^i}
 -\rho^i_b\pd{\psi^\alpha_c}{x^i}
 +\psi^\alpha_aC^a_{bc}\right) 
\bar{e}^b\wedge\bar{e}^c.
\end{align*}
Finally, for $\theta=\V^\beta_b=d y^\beta_b$ we also get an identity 
$$
d\Psi\pb\V^\alpha_a-\Psi\pb d\V^\alpha_a
=d\Psi\pb dy^\alpha_a
=dd\Psi\pb y^\alpha_a=0,
$$
provided that the admissibility conditions hold.

Therefore, $\Psi$ is a morphism if and only if in addition to (a) and (b) it satisfies the equation
\begin{enumerate}
\item[(c)]
$\p{\psi}^\alpha_{b|c}-\p{\psi}^\alpha_{c|b}
-C^\alpha_{\beta\gamma}\psi^\beta_b\psi^\gamma_c
 -C^\alpha_{b\gamma}\psi^\gamma_c
 +C^\alpha_{c\beta}\psi^\beta_b
 +\psi^\alpha_aC^a_{bc}
 +C^\alpha_{bc}=0$.
\end{enumerate}

Let us consider a semiholonomic section $\check{\Psi}$, so that $\psi^\alpha_a=\phi^\alpha_a$. Then the above condition (a) reads $\p{\phi}^A_{|a}=\rho^A_a+\rho^A_\alpha\phi^\alpha_a$ which is but the admissibility conditions for the map $\Phi$, and condition (b) reads $\p{\phi}^\beta_{b|a}=\psi^\beta_{ba}$, which is just expressing that $\Psi=\Phi\spi$. This proves the first statement.

Moreover, the section $\Psi$ is a morphism if it is admissible, and therefore it is of the form $\Psi=\Phi\spi$ for $\Phi$ admissible, and in addition it satisfies condition (c) which is just expressing that the admissible map $\Phi$ is a morphism. This proves the second statement.
\end{proof}

Obviously the relation between $\Psi$ and $\Phi$ is by projection (both projections): $\prol{\prJM}\circ\Psi(n)=\Phi_n=(\pi_1)_{10}(\Psi(n))$ for all $n\in N$.

\begin{corollary}
Let $\Phi$ be an admissible section of $\pi$. Then $\Phi$ is a morphism of Lie algebroids if and only if $\check{\Phi}\spi$ is holonomic
\end{corollary}
\begin{proof}
Notice first that, by the above proposition, if $\Phi$ is admissible then $\Psi=\Phi\spi$ is semiholonomic and admissible, i.e. satisfies conditions (a) y (b) in the proof of the theorem. The result follows by noticing that $\Phi\spi$ is holonomic if and only if $\mathcal{M}^\alpha_{ab}=0$, which is just condition (c).
\end{proof}



\section{The Lagrangian formalism}
\label{Lagrangian}

We consider a Lagrangian function $L\in\cinfty{\Jpi}$ and a fixed nowhere-vanishing form $\omega\in\ext[r]{F}$ of maximal degree $r=\operatorname{rank}(F)$, to which we refer as the volume form. We will denote by $\tilde{\omega}=\pi_1\pb\omega$ the pullback of the volume form to $\prol[E]{\Jpi}$ by the projection $\pi_1$. The product $\L=L\tilde{\omega}\in\ext[r]{\prol[E]{\Jpi}}$ will be called the Lagrangian density. 

We will define in this section the Cartan forms associated to the Lagrangian and, in terms of them, we will get a system of partial differential equations which will be called the Euler-Lagrange partial differential equations. In order to do that we will consider the analogs of the the vertical lifting and the vertical endomorphism in the standard theory.

\subsection*{Vertical lifting}
The bundle $\Jpi$ being affine has a well defined vertical lifting map. If $\phi\in\Jpi[m]$ then we have a map $\map{{\,\,}^V_\phi}{\Vpi[m]}{T_\phi\Jpi}$,  given by
$$
\psi^V_\phi f=\frac{d}{dt}f(\phi+t\psi)\Big|_{t=0},
$$
for every function $f\in\cinfty{\Jpi}$. Moreover, in the special case of $\Jpi$, we have some additional structure. Indeed, an element $\phi$ is but a splitting of an exact sequence, and we have the associated vertical projector $\map{v_\phi}{E_m}{K_m\subset E_m}$ given by  $v_\phi(a)=a-\phi(\prEF(a))$. This allows to extend the vertical lifting operation to elements $\varphi\in\Lpi[m]$, by  $(\v_\phi\circ\varphi)^V_\phi$, which is well defined because $\textrm{Im}(v_\phi)\subset K_{m}$. 

\begin{definition}
The map $\map{\xi^V}{\prJM^*(\Lpi)}{\prol[E]{\Jpi}}$ given by
$$
\xi^V(\phi,\varphi)=\bigl(\phi,0,(v_\phi\circ\varphi)^V_\phi\bigr).
$$
is said to be the  \emph{vertical lifting} map.
\end{definition}

In coordinates, if $\phi=(e_a+y^\alpha_ae_\alpha)\otimes\bar{e}^a$ and $\varphi=(\varphi^b_a e_b +\varphi^\alpha_ae_\alpha)\otimes\bar{e}^a$ then $v_\phi=e_\alpha\otimes(e^\alpha-y^\alpha_ae^a)$
and $\xi\spV(\phi,\varphi)= (\varphi^\alpha_a-y^\alpha_b\varphi^b_a)\V_\alpha^a$.

This construction will be generally used in the case of a map $\varphi$ of the form $\varphi=b\otimes\lambda$, with $b\in E_m$ and $\lambda\in F_n$. In such case
$$
\xi^V(\phi,b\otimes\lambda)=\bigl(\phi,0,(v_\phi(b)\otimes\lambda)^V_\phi\bigr),
$$
and its coordinate expression is $\xi\spV(\phi,b\otimes\lambda)= \lambda_a(b^\alpha-y^\alpha_cb^c)\V_\alpha^a$.

\subsection*{Vertical endomorphism}
Given a section $\nu\in\sec{F^*}$ we can define a linear map $\map{S_\nu}{\prol[E]{\Jpi}}{\prol[E]{\Jpi}}$, known as the \emph{vertical endomorphism}, by means of projection and vertical lifting. If $Z=(\phi,a,V)\in\prol[E]{\Jpi}$ then
$$
S_\nu(Z)
=\xi^V\bigl(\phi,a\otimes\nu_n\bigr),
$$
where $m=\prJM(\phi)$ and $n=\prMN(m)$.

From this expression it is clear that $S_\nu$ depends linearly in $\nu$, a fact which helps to find the coordinate expression of $S_\nu$. If $\nu=\nu_a\bar{e}^a$ and we denote $S^a\equiv S_{\bar{e}^a}$ then $S_\nu=\nu_a S^a$ and
\begin{align*}
S^a(\X_b(\phi))
&=S^a\Bigl(\phi,e_b,\rho^i_b\pd{}{x^i}+\rho^A_b\pd{}{u^A}\Bigr)
=\xi^V(\phi,v_\phi(e_b)\otimes\bar{e}^a)\\
&=-y^\alpha_b\xi^V(\phi,e_\alpha\otimes\bar{e}^a)
=-y^\alpha_b\V^a_\alpha(\phi),
\end{align*}
$$
S^a(\X_\beta(\phi))
=S^a\Bigl(\phi,e_\beta,\rho^A_\beta\pd{}{u^A}\Bigr)
=\xi^V(\phi,v_\phi(e_\beta)\otimes\bar{e}^a)
=\V^a_\beta(\phi)
$$
and
$$
S^a(\V^b_\beta(\phi))
=S^a\biggl(\phi,0,\pd{}{y^\beta_b}\biggr)
=\xi^V(\phi,0)
=0.
$$
Thus the coordinate expression of $S^a$ is $S^a=(\X^\alpha-y^\alpha_b\X^b)\otimes\V^a_\alpha =\theta^\alpha\otimes\V^a_\alpha$, where we recall that $\theta^\alpha=\X^\alpha-y^\alpha_b\X^b$ are the elements of a basis of contact 1-forms.
Therefore the expression of the vertical endomorphism is 
$$
S=\theta^\alpha\otimes\bar{e}_a\otimes\V^a_\alpha.
$$
As usual we will denote $S_\omega$ the $r$-form obtained by contraction (of the second tensorial factor) with the volume form $\omega$. In coordinates
$$
S_\omega=\theta^\alpha\wedge\omega_a\otimes\V^a_\alpha.
$$

\subsection*{Cartan forms}
With the help of the vertical endomorphism we define the Lagrangian \emph{multimomentum $r$-form} $\Theta_L\in\ext[r]{\prol[E]{\Jpi}}$ by
$$
\Theta_L=S_\omega(dL)+L\omega,
$$
which in local coordinates reads
$$
\Theta_L=\pd{L}{y^\alpha_a}\theta^\alpha\wedge\omega_a+L\omega.
$$
By taking the differential of the multimomentum form we have the Lagrangian \emph{multisymplectic $(r+1)$-form} $\Omega_L=-d\Theta_L$. In coordinates we have
\begin{align*}
\Omega_L&=\frac{1}{2}\left[
 \pd{^2L}{y^\beta_a\partial u^B}\rho^B_\gamma
-\pd{^2L}{y^\gamma_a\partial u^B}\rho^B_\beta
+\pd{L}{y^\alpha}C^\alpha_{\beta\gamma}
\right]\theta^\beta\wedge\theta^\gamma\wedge\omega_a
+\pd{^2L}{y^\alpha_a\partial y^\beta_b}
  \theta^\alpha\wedge\V^\beta_b\wedge\omega_a+\\
&{}+\left[
 \pd{^2L}{y^\alpha_a\partial u^B}(\rho^B_a+\rho^B_\beta y^\beta_a)
+\pd{^2L}{y^\alpha_a\partial x^i}\rho^i_a
+\pd{L}{y^\alpha_a}C^b_{ba}
-\pd{L}{y^\gamma_a}Z^\gamma_{a\alpha}
-\pd{L}{u^A}\rho^A_\alpha
\right]\theta^\alpha\wedge\omega.
\end{align*}

\subsection*{Euler-Lagrange equations}
By a solution of the Hamiltonian system defined by a Lagrangian $L$ we mean a morphism $\map{\Phi}{F}{E}$ such that its prolongation satisfies the equation
$$
\Phi\spi{}\pb(i_X\Omega_L)=0,
$$
for every section $X$ of $\prol[E]{\Jpi}$ vertical over $F$. The above equations are said to be the \emph{Euler-Lagrange equations} for the Lagrangian $L$. We will obtain in this section the local expression of this equations and the system of partial differential equations which defines. 

We will start from a slightly more general point of view and we will look for sections $\Psi$ of $\pi_1$ such that they satisfy the so called \emph{De Donder equations},
$$
\Psi\pb(i_X\Omega_L)=0,
$$
for every section $X$ of $\prol[E]{\Jpi}$ vertical over $F$. Later on we will impose that $\Psi$ is admissible and finally a morphism.

Since $\ul{\Psi}$ is a section of $\ul{\pi_1}$ it has an associated section $\Phi\in\Sec{\pi}$, such that $\check{\ol{\Phi}}=\ul{\Psi}$. In coordinates $\Psi$ is of the form $\Psi=(\X_a+\psi^\alpha_a\X_\alpha+\psi^\beta_{ba}\V^b_\beta)\otimes\bar{e}^a$ and  the base map $\Phi$ is of the form $\Phi(x)=(x^i,\phi^A(x),\phi^\alpha_a(x))$.
Therefore we have
\begin{align*}
&\Psi\pb x^i=x^i
&&\Psi\pb\X^a=\bar {e}^a\\
&\Psi\pb u^A=\phi^A
&&\Psi\pb\X^\beta=\psi^\beta_a\bar {e}^a\\
&\Psi\pb y^\alpha_a=\phi^\alpha_a
&&\Psi\pb\V^\beta_b=\psi^\beta_{ba}\bar {e}^a.
\end{align*}

Taking $X=\V^b_\beta$ in the De Donder equations we get
$$
i_X\Omega_L=-\pd{^2L}{y^\alpha_a\partial y^\beta_b}\theta^\alpha\wedge\omega_a
$$
and thus
$$
\Psi\pb(i_X\Omega_L)=
-\pd{^2L}{y^\alpha_a\partial y^\beta_b}(\psi^\alpha_a-\phi^\alpha_a)\wedge\omega
$$
which vanishes if and only if
$$
\pd{^2L}{y^\alpha_a\partial y^\beta_b}(\psi^\alpha_a-\phi^\alpha_a)=0.
$$

We will say that the Lagrangian is \emph{regular} if the linear map $\xi^\alpha_a\mapsto\pd{^2L}{y^\alpha_a\partial y^\beta_b}\xi^\alpha_a$ is regular at every point. In the case of a regular Lagrangian the solution of the above equation is just $\psi^\alpha_a=\phi^\alpha_a$, that is $\Psi$ is semiholonomic. On the other hand, if the section $\Psi$ is semiholonomic, the above conditions are always satisfied, that is 

\begin{proposition}\label{semiholonomic.and.vertical}
If $\Psi$ is a semiholonomic section then $\Psi\pb(i_X\Omega_L)=0$ for every  section $X$ of $\prol[E]{\Jpi}$ vertical over $E$.
\end{proposition}

In general we will be treating with singular Lagrangians. Therefore, in what follows we will consider only semiholonomic sections. 

\medskip

Taking now $X=\X_\alpha$ we obtain
\begin{align*}
i_X\Omega_L=
&=\frac{1}{2}\left[
 \pd{^2L}{y^\alpha_a\partial u^B}\rho^B_\beta
-\pd{^2L}{y^\beta_a\partial u^B}\rho^B_\alpha
+\pd{L}{y^\gamma}C^\gamma_{\alpha\beta}
\right]\theta^\beta\wedge\omega_a+
\pd{^2L}{y^\alpha_a\partial y^\beta_b}\V^\beta_b\wedge\omega_a+\\
&{}+\left[
 \pd{^2L}{y^\alpha_a\partial u^B}(\rho^B_a+\rho^B_\beta \Phi^\beta_a)
+\pd{^2L}{y^\alpha_a\partial x^i}\rho^i_a
+\pd{L}{y^\alpha_a}C^b_{ba}
-\pd{L}{y^\gamma_a}Z^\gamma_{a\alpha}
-\pd{L}{u^A}\rho^A_\alpha
\right]\omega.
\end{align*}
and thus 
\begin{align*}
\Psi\pb(i_X\Omega_L)&=
\Biggl[
\pd{^2L}{y^\alpha_a\partial y^\beta_b}\psi^\beta_{ba}+
 \pd{^2L}{y^\alpha_a\partial u^B}(\rho^B_a+\rho^B_\beta\phi^\beta_a)
+\pd{^2L}{y^\alpha_a\partial x^i}\rho^i_a+\hbox to 2cm{\hss}\\
&\hbox to 5cm{\hss}{}+\pd{L}{y^\alpha_a}C^b_{ba}
-\pd{L}{y^\gamma_a}Z^\gamma_{a\alpha}
-\pd{L}{u^A}\rho^A_\alpha
\Biggr]\omega,
\end{align*}
so that the coefficients $\psi^\alpha_{ba}$ must satisfy
$$
\pd{^2L}{y^\alpha_a\partial y^\beta_b}\psi^\beta_{ba}+
 \pd{^2L}{y^\alpha_a\partial u^B}(\rho^B_a+\rho^B_\beta\phi^\beta_a)
+\pd{^2L}{y^\alpha_a\partial x^i}\rho^i_a
+\pd{L}{y^\alpha_a}C^b_{ba}
-\pd{L}{y^\gamma_a}Z^\gamma_{a\alpha}
-\pd{L}{u^A}\rho^A_\alpha=0.
$$

In addition to this equations we have the admissibility condition. As a consequence of proposition~\ref{admissible-prolongation}, if $\Psi$ is semiholonomic and admissible, then $\Psi=\Phi\spi$ with $\Phi$ admissible. Thus
$$
\p{\phi}^A_{|a}=\rho^A_a+\rho^A_\alpha \phi^\alpha_a
\qquand
\p{\phi}^\beta_{b|a}=\psi^\beta_{ba},
$$
so that the first three terms in the above equation are but the expression of the total derivative of $\partial L/\partial y^\alpha_a$, so that we finally get 
$$
\left(\pd{L}{y^\alpha_a}\right)'_{|a}
+\pd{L}{y^\alpha_a}C^b_{ba}
-\pd{L}{y^\gamma_a}Z^\gamma_{a\alpha}
-\pd{L}{u^A}\rho^A_\alpha=0.
$$

If we finally impose that $\Psi$ is not only admissible but a morphism (i.e. $\Phi$ is a morphism) we get the system of partial differential equations
\begin{align*}
&\rho^i_a\pd{u^A}{x^i}=\rho^A_a+\rho^A_\alpha y^\alpha_a\qquad\qquad
\rho^i_a\pd{y^\beta_b}{x^i}=y^\beta_{ba}\\
&\bigr(y^\gamma_{ab}+C^\gamma_{b\alpha}y^\alpha_a\bigr)
-\bigr(y^\gamma_{ba}+C^\gamma_{a\beta}y^\beta_b\bigr)
 -C^\gamma_{\alpha\beta}y^\alpha_ay^\beta_b+y^\gamma_c C^c_{ab}
+C^\gamma_{ab}=0\\
&\left(\pd{L}{y^\alpha_a}\right)'_{|a}
+\pd{L}{y^\alpha_a}C^b_{ba}
-\pd{L}{y^\gamma_a}Z^\gamma_{a\alpha}
-\pd{L}{u^A}\rho^A_\alpha=0,
\end{align*}
which are said to be the Euler-Lagrange partial differential equations.

\subsection*{Alternative expression of $\Omega_L$ and the Euler-Lagrange form}
If we introduce the forms $\theta^\beta_b=\V^\beta_b-y^\beta_{ba}\X^a$, which are forms on $\prol[E]{\Jtpi}$, then the pullback of $\Omega_L$ to $\prol[E]{\Jtpi}$ is (we omit the pullback in the notation)
\begin{align*}
\Omega_L&=\frac{1}{2}\left[
 \pd{^2L}{y^\beta_a\partial u^B}\rho^B_\gamma
-\pd{^2L}{y^\gamma_a\partial u^B}\rho^B_\beta
+\pd{L}{y^\alpha}C^\alpha_{\beta\gamma}
\right]\theta^\beta\wedge\theta^\gamma\wedge\omega_a
+\pd{^2L}{y^\alpha_a\partial y^\beta_b}
  \theta^\alpha\wedge\theta^\beta_b\wedge\omega_a+\\
&{}+\left[
 \pd{^2L}{y^\alpha_a\partial y^\beta_b}y^\beta_{ba}
+\pd{^2L}{y^\alpha_a\partial u^B}(\rho^B_a+\rho^B_\beta y^\beta_a)
+\pd{^2L}{y^\alpha_a\partial x^i}\rho^i_a
+\pd{L}{y^\alpha_a}C^b_{ba}
-\pd{L}{y^\gamma_a}Z^\gamma_{a\alpha}
-\pd{L}{u^A}\rho^A_\alpha
\right]\theta^\alpha\wedge\omega.
\end{align*}
and the last term is the Euler-Lagrange $(r+1)$-form
$$
\delta L=\left[
 \left(\pd{L}{y^\alpha_a}\right)'_{|a}
+C^b_{ba}\left(\pd{L}{y^\alpha_a}\right)
-\pd{L}{y^\gamma_a}Z^\gamma_{a\alpha}
-\pd{L}{u^A}\rho^A_\alpha
\right]\theta^\alpha\wedge\omega.
$$
The Euler-Lagrange equations can be written as $\delta L=0$, which define a subset (a submanifold under some regularity conditions) of $\Jtpi$ on which the jet prolongation $\Phi\spi$ of a solution $\Phi$ has to take value.


\section{ Variational Calculus}
\label{variational}
In this section we will show that, in the case $F=TN$, our formalism  admits a variational description, i.e. the Euler-Lagrange equations obtained by the multisymplectic formalism are the equations for critical sections of a constrained variational problem~\cite{DGA04}.

\subsection*{Variational problem}
We consider the following variational problem: Given a Lagrangian function $L\in\cinfty{\Jpi}$ and a volume form  $\omega$ on $N$, find the critical points of the action functional 
$$
\CMcal{S}(\Phi)=\int_N (L\circ\check{\ol{\Phi}})\,\omega
$$
defined on the set of morphisms sections of $\pi$, that is, on the set $\Mor{\pi}$.

The above variational problem is a constrained problem, not only because the condition $\pi\circ\Phi=\id_F$, which can be easily solved, but because we are restricting $\Phi$  to be a morphism of Lie algebroids, which is a condition on the derivatives of $\Phi$. We will explicitly find some curves in the space of morphisms which will allow us to get the equations satisfied by the critical sections.

In order to use all our geometric machinery, we will reformulate the variational problem. If we consider the Lagrangian density $\L=L\tilde{\omega}\in\ext{\prol[E]{\Jpi}}$, then we have that 
$$
(L\circ\check{\ol{\Phi}})\,\omega=\Phi\spi{}\pb\L.
$$
Indeed, the base map of $\Phi\spi$ is $\check{\ol\Phi}$ so that we have $\Phi\spi{}\pb L=L\circ\check{\ol{\Phi}}$. Moreover, since $\Phi\spi$ is a section of $\pi_1$, we have that $\Phi\spi{}\pb\tilde{\omega}=\omega$. Therefore we can write the action as $\CMcal{S}(\Phi)=\int_N\Phi\spi{}\pb\L$. 

On the other hand, since the difference between $\L$ and $\Theta_L$ is a contact form and $\Phi\spi$ is semiholonomic, we have that $\Phi\spi{}\pb\L=\Phi\spi{}\pb\Theta_L$. Therefore the action functional can be rewritten as
$$
\CMcal{S}(\Phi)=\int_N\Phi\spi{}\pb\Theta_L.
$$

In order to find admissible variations we consider sections of $E$ and the associated flow. With the help of this flow we can transform morphisms of Lie algebroids into morphisms of Lie algebroids as it is explained in the next subsection.

\subsection*{Jet prolongation of maps and sections}
In this subsection we return momentarily to the general case of a general Lie algebroid $F$. 

Consider a vector bundle map $\Psi=(\ol{\Psi},\ul{\Psi})$ from $E$ to $E$ which induces the identity in $F$, that is $\pi\circ\Psi=\pi$, or explicitly $\prEF\circ\ol{\Psi}=\prEF$ and $\prMN\circ\ul{\Psi}=\prMN$. The map $\ol\Psi$ induces a map between jets $\map{\check{\ol{\Psi}}}{\Jpi}{\Jpi}$ by composition, $\check{\ol{\Psi}}(\phi)=\ol{\Psi}\circ\phi$. The map $\check{\ol{\Psi}}$ is well defined since $\prEF\circ(\check{\ol{\Psi}}(\phi))= \prEF\circ\ol\Psi\circ\phi = \prEF\circ\phi=\id$, that is, the image of $\check{\ol{\Psi}}$ is in $\Jpi$. Moreover, $\check{\ol{\Psi}}$ is a fiberwise affine map over the map $\ul{\Psi}$. Thus we have the affine bundle map $\check{\Psi}=(\check{\ol{\Psi}},\ul{\Psi})$ from $\map{\ul{\pi_{10}}}{\Jpi}{M}$ to itself. By tangent prolongation we have the bundle map $\Psi\spi=\prol[\Psi]{\check{\Psi}}$ from $\map{\tau^E_{\Jpi}}{\prol[E]{\Jpi}}{\Jpi}$ to itself  which will be called the \emph{jet-prolongation} of $\Psi$.

\begin{proposition}
Let $\Phi=(\ol{\Phi},\ul{\Phi})$ be a section of $\pi$ and $\Psi$ a vector bundle map over the identity in $F$. We consider the transformed map $\Phi'=\Psi\circ\Phi$, i.e. $\ol{\Psi'}=\ol{\Psi}\circ\ol{\Phi}$ and $\ul{\Phi'}=\ul{\Psi}\circ\ul{\Phi}$. Then 
\begin{enumerate}
\item The maps $\check{\Phi}$ and $\check{\Phi}'$ associated to $\Phi$ and $\Phi'$ are related by  $\check{\Phi}'=\check{\Psi}\circ\check{\Phi}$. 
\item The jet prolongation of $\Phi$ and of $\Phi'$ are related by ${\Phi'}\spi=\Psi\spi\circ\Phi\spi$.
\end{enumerate}
\end{proposition}
\begin{proof}
Indeed, the base maps are $\ul{\Phi}'=\ul{\Psi}\circ\ul{\Phi}$, and  for every $n\in N$,
$$
\check{\ol{\Phi}}'(n)=\ol{\Phi}'_n =\ol{\Psi}\circ\ol\Phi_n=\ol{\Psi}\circ(\check{\ol{\Phi}}(n))
=\check{\ol{\Psi}}(\check{\ol{\Phi}}(n)),
$$
which proves the first. For the second we have
$$
{\Phi'}\spi
=\prol[\Phi']{\check{\Phi}'}
=\prol[\Psi\circ\Phi]{\check{\Psi}\circ\check{\Phi}}
=\prol[\Psi]{\check{\Psi}}\circ\prol[\Phi]{\check{\Phi}}
=\Psi\spi\circ{\Phi}\spi,
$$
where we have used the composition property for $E$-tangent prolongations.
\end{proof}

We will say that a section $\sigma$ of $E$ is $\pi$-vertical if it projects to the zero section on $F$, that is $\prEF\circ\sigma=0$. Let $\sigma$ be a $\pi$-vertical section of $E$, and consider the flow $\Psi_s=(\ol\Psi_s,\ul\Psi_s)$ of $\sigma$. We recall that  $\ul\Psi_s$ is the flow of the vector field $\rho(\sigma)$. 

The associated map $\check{\ol\Psi}_s$ is a flow on the manifold $\Jpi$ and we denote by  $\check{X}_\sigma\in\vectorfields{\Jpi}$ its infinitesimal generator. It follows that $T\prJM\circ\check{X}_\sigma=\rho(\sigma)\circ\prJM$ and thus for every $m\in M$ and every $\phi\in\Jpi[m]$ we have that $(\phi,\sigma(m),\check{X}_\sigma(\phi))$ is an element of $\prol[E]{\Jpi}$ at the point $\phi$. Therefore we have defined a section $\sigma\spi\in\sec{\prol[E]{\Jpi}}$ given by 
$$
\sigma\spi(\phi)=(\phi,\sigma(m),\check{X}_\sigma(\phi)),
$$
where $m=\prJM(\phi)$. This section is said to be the \emph{jet-prolongation} of $\sigma\in\sec{E}$ or the \emph{complete lift} to $\Jpi$ of the section $\sigma$. 

\begin{proposition}\label{complete.lift}
The first jet prolongation of a $\pi$-vertical section satisfies the following properties 
\begin{enumerate}
\item\label{a}The jet prolongation $\sigma\spi$ of $\sigma$ is projectable and projects to $\sigma$.
\item\label{b}The jet prolongation $\sigma\spi$ of $\sigma$ preserves the contact ideal, that is, if $\theta$ is a contact form then $d_{\sigma\spi}\theta$ is a contact form.
\item\label{c}The jet prolongation $\sigma\spi$ of $\sigma$ is the only section of $\prol[E]{\Jpi}$ which satisfies the above two properties~(\ref{a}) and~(\ref{b}).
\item\label{d}The flow of the jet prolongation $\sigma\spi$ is $\Psi_s\spi$, where $\Psi_s$ is the flow of the section $\sigma$. 
\item\label{e}If the local expression of the section $\sigma$ is $\sigma=\sigma^\alpha e_\alpha$, then the local expression of the the jet prolongation of $\sigma$ is 
$$
\sigma\spi=
\sigma^\alpha\X_\alpha
+(\p{\sigma}^\alpha_{|a}+Z^\alpha_{a\beta}\sigma^\beta)\V^a_\alpha
$$
\end{enumerate}
\end{proposition}
\begin{proof}
Property~(\ref{a}) is obvious from the definition of $\sigma\spi$ and  property~(\ref{d}) follows from proposition~\ref{flow.projectable}. 

To prove~(\ref{b}) we notice that if $\Psi$ is a bundle map and if $\beta$ is a section of $E^*$ and $\breve{\beta}$ is the associated contact form then $\Psi\spi{}\pb\breve{\beta}=\breve{\Psi\pb\beta}$, the contact form associated to $\Psi\pb\beta$. Then~(\ref{b}) follows since the flow of $\sigma\spi$ is $\Psi_s\spi$. 

In local coordinates, from property~(\ref{a}) we have that $\sigma\spi$ has the form $\sigma\spi=\sigma^\alpha\X_\alpha+\sigma^\alpha_a\V^a_\alpha$, for some functions $\sigma^\alpha_a$. Then taking the Lie derivative of a contact form $\theta^\alpha=\X^\alpha-y^\alpha_a\X^a$ we get
\begin{align*}
d_{\sigma\spi}\theta^\alpha
&=i_{\sigma\spi}d\theta^\alpha+di_{\sigma\spi}\theta^\alpha \\
&=C^\alpha_{\beta\gamma}\sigma^\gamma\theta^\beta+Z^\alpha_{a\gamma}\sigma^\gamma\X^a-\sigma^\alpha_a\X^a+d\sigma^\alpha\\
&=(Z^\alpha_{a\gamma}\sigma^\gamma+\p{\sigma}^\alpha_{|a}-\sigma^\alpha_a)\X^a+\left(\rho^A_\beta\pd{\sigma^\alpha}{u^A}+C^\alpha_{\beta\gamma}\sigma^\gamma\right)\theta^\beta.
\end{align*}
which is a contact form if and only if $\sigma^\alpha_a=\p{\sigma}^\alpha_{|a}+Z^\alpha_{a\gamma}\sigma^\gamma$. This proves~(\ref{e}). Finally, property~(\ref{c}) follows from the above coordinate expression.
\end{proof}

\begin{proposition}
If $f\in\cinfty{N}$ then $(f\sigma)\spi=f\sigma\spi+(df\otimes\sigma)\spV$.
\end{proposition}
\begin{proof}
Both $(f\sigma)\spi$ and $f\sigma\spi$ project to $f\sigma$, so that they differ in a vertical section $V=(f\sigma)\spi-f\sigma\spi$. From the coordinate expression of $\sigma\spi$ we have that $V=\p{f}_{|a}\sigma^\alpha\V_\alpha^a$ which is but the local expression of $(df\otimes\sigma)\spV$.
\end{proof}

\begin{remark}
In more generality, one can define the complete lift of any  section $\sigma$ of $E$ which is projectable to a section $\bar{\sigma}$ of $F$, not necessarily the zero section. Indeed, we just need to define the flow of the vectorfield $\rho^1(\sigma\spi)$. The flow $\Psi_s$ of the section $\sigma$ projects to the flow $\Phi_s$ of the section $\bar{\sigma}$, that is $\pi\circ\Psi_s=\Phi_s\circ\pi$. Then we can define the maps $\map{\mathcal{J}\Psi_s}{\Jpi}{\Jpi}$ by means of $\mathcal{J}\Psi_s(\phi)=\ol{\Psi}_s\circ\phi\circ\ol{\Phi}_{-s}$. This maps define a local flow on the manifold $\Jpi$ whose infinitesimal generator is (by definition) the vectorfield $\rho^1(\sigma\spi)$.

Even more generally, one can define the complete lift of a nonprojectable section $\sigma$ of $E$ by imposing the first two conditions in proposition~\ref{complete.lift}. We will not study such constructions in this paper since we will not need it. 
\end{remark}

\subsection*{Equations for critical sections}
Going back to our variational problem for $F=TN$, let $\Phi$ be a critical point of $\CMcal{S}$. An admissible variation of $\Phi$ is but a curve $\mathcal{M}(\pi)$ starting at $\Phi$, that is a map $s\mapsto \Phi_s$ such that $\Phi_s$ is (for every fixed $s$) a section of $\pi$ and a morphism of Lie algebroids. To find one of such curves we consider a $\pi$-vertical section $\sigma$ of $E$ (thus its flow $\map{\Psi_s}{E}{E}$ projects to the identity in $F=TN$) with compact support. It follows that, for every fixed $s$, the bundle map $\Phi_s=\Psi_s\circ\Phi$ is a section of $\pi$ and a morphism of Lie algebroids, that is, $s\mapsto\Phi_s$ is a curve in $\mathcal{M}(\pi)$. Thus we have that
\begin{align*}
0&=\frac{d}{ds}\CMcal{S}(\Phi_s)\at{s=0}
=\int_N\frac{d}{ds}\Phi_s\spi{}\pb\Theta_L\at{s=0}
=\int_N\frac{d}{ds}{}\bigl(\Psi_s\spi\circ\Phi\spi)\pb\Theta_L\bigr)\at{s=0}\\
&=\int_N\frac{d}{ds}\Phi\spi{}\pb\bigl(\Psi_s\spi{}\pb\Theta_L\bigr)\at{s=0}
=\int_N\Phi\spi{}\pb\frac{d}{ds}\bigl(\Psi_s\spi{}\pb\Theta_L\bigr)\at{s=0}
=\int_N\Phi\spi{}\pb d_{\sigma\spi}\Theta_L,
\end{align*}
where we have used that $\Phi_s\spi=\Psi_s\spi\circ\Phi\spi$ and that $\Psi_s\spi$ is the flow of the jet prolongation $\sigma\spi$ of the section $\sigma$.

Using that $d_{\sigma\spi}=i_{\sigma\spi}\circ d+d\circ i_{\sigma\spi}$ we have that $d_{\sigma\spi}\Theta_L=i_{\sigma\spi}d\Theta_L+di_{\sigma\spi}\Theta_L =-i_{\sigma\spi}\Omega_L+d i_{\sigma\spi}\Theta_L$, so that 
$$
0=\int_N\Phi\spi{}\pb d_{\sigma\spi}\Theta_L
=-\int_N\Phi\spi{}\pb i_{\sigma\spi}\Omega_L+\int_N\Phi\spi{}\pb d i_{\sigma\spi}\Theta_L.
$$
The second term vanishes; indeed $\Theta_L$ is semibasic, so that $i_{\sigma\spi}\Theta_L$ depends only on the values of $\sigma$, not on their derivatives, and since $\sigma$ has compact support we have that 
$$
\int_N\Phi\spi{}\pb d i_{\sigma\spi}\Theta_L
=\int_N d[\Phi\spi{}\pb i_{\sigma\spi}\Theta_L]=0.
$$
Therefore we get
$$
0=\int_N\Phi\spi{}\pb i_{\sigma\spi}\Omega_L.
$$
We now prove that this implies $\Phi\spi{}\pb i_{\sigma\spi}\Omega_L=0$ for every $\sigma$. Indeed, if we take the section $f\sigma$, for $f\in\cinfty{N}$, then $(f\sigma)\spi=f\sigma\spi+(df\otimes\sigma)^V$, so that 
$$
\Phi\spi{}\pb i_{(f\sigma)\spi}\Omega_L=
f\Phi\spi{}\pb i_{\sigma\spi}\Omega_L+\Phi\spi{}\pb i_{(df\otimes\sigma)^V}\Omega_L.
$$
But $\Phi\spi$ is semiholonomic and $(df\otimes\sigma)^V$ is vertical, so that by proposition~\ref{semiholonomic.and.vertical} we have $\Phi\spi{}\pb i_{(df\otimes\sigma)^V}\Omega_L=0$. Therefore
$$
0=\int_N\Phi\spi{}\pb i_{(f\sigma)\spi}\Omega_L
=\int_Nf[\Phi\spi{}\pb i_{\sigma\spi}\Omega_L],
$$
for every function $f\in\cinfty{N}$. From the fundamental theorem of the Calculus of Variations we get that $\Phi\spi{}\pb i_{\sigma\spi}\Omega_L=0$ for every $\sigma$. 

Finally we show that this condition is equivalent to $\Phi\spi{}\pb i_X\Omega_L=0$ for every $\pi_1$-vertical section $X$ of $\prol[E]{\Jpi}$. Indeed, the above equation is tensorial in $X$, so that it is equivalent to the same condition for projectable sections $X$. But a if $X$ is projectable and projects to $\sigma$, then $X=\sigma\spi+V$ for some vertical section. Thus 
$\Phi\spi{}\pb i_X\Omega_L=\Phi\spi{}\pb i_{\sigma\spi}\Omega_L+\Phi\spi{}\pb i_V\Omega_L$, but the last term vanishes because $V$ is vertical and $\Phi\spi$ is semiholonomic.

Therefore, in the case of a tangent bundle $F=TN$ with the standard Lie algebroid structure, we have proved that 
\begin{theorem}
The following conditions are equivalent:
\begin{enumerate}
\item A morphism $\Phi$ is a critical section of $\CMcal{S}$.
\item In local coordinates such that the volume form is $\omega=dx^1\wedge\cdots\wedge dx^r$ and $\rho^i_a=\delta^i_a$, the components $y^\alpha_a$ of a vector bundle map $\Phi\in\sec{\pi}$  satisfy the system of partial differential equations
\begin{align*}
&\pd{u^A}{x^a}=\rho^A_a+\rho^A_\alpha y^\alpha_a\\
& \pd{y^\alpha_a}{x^b}-\pd{y^\alpha_b}{x^a}
 +C^\alpha_{b\gamma}y^\gamma_a-C^\alpha_{a\gamma}y^\gamma_b
 +C^\alpha_{\beta\gamma}y^\beta_by^\gamma_a
 +C^\alpha_{ab}=0\\
&\frac{d\,\,}{dx^a}\left(\pd{L}{y^\alpha_a}\right)
=\pd{L}{y^\gamma_a}Z^\gamma_{a\alpha}
+\pd{L}{u^A}\rho^A_\alpha.
\end{align*}
\item A morphism $\Phi$ satisfies the Euler-Lagrange equations $\Phi\spi{}\pb i_X\Omega_L=0$, for every $\pi_1$-vertical section $X$.
\end{enumerate}
\end{theorem}
\begin{proof}
The equivalence of (1) and (3) has been already established. In order to prove that (2) is equivalent to (3) we just have to notice that 
$$
C^a_{bc}=0
\qquand
\p{f}_{|a}=\frac{df}{dx^a}
$$
because $e_i=\partial/\partial x^i$ is a coordinate basis.
\end{proof}

\subsection*{Noether's theorem}
Noether's theorem is a consequence of the existence of a variational description of the problem. In the standard case~\cite{Barna-L}, when the Lagrangian is invariant by the first jet prolongation of a vertical vectorfield $Z$ then the Noether current $J=i_{Z\spi}\Theta_L$ is a conserved current. By a conserved current we mean a $(r-1)$-form such that its pullback by any solution of the Euler-Lagrange equations is a closed form on the base manifold $N$. Therefore the integral of the $(r-1)$-form $J$ over any closed $(r-1)$-dimensional submanifold vanishes.

We will show that a similar statement can be obtained for a field theory over general Lie algebroids $F$ and $E$: for every symmetry of the Lagrangian we get a form which is closed over every solution of the Euler-Lagrange equations. But now the word \textsl{conserved} has only a partial meaning, because when the Lie algebroid $F$ is not $TN$ we do not have a volume integral neither a Stokes theorem.

In what follows in this subsection $F$ is a general Lie algebroid, not necessarily a tangent bundle.

Following the steps in the derivation of the equation for critical sections, if $\sigma$ is any $\pi$-vertical section of $E$ (we do not need it to be compactly supported) we have that 
$$
\Phi\spi{}\pb d_{\sigma\spi}\L 
= d[\Phi\spi{}\pb i_{\sigma\spi}\Theta_L]-\Phi\spi{}\pb i_{\sigma\spi}\Omega_L,
$$
It follows that if $\Phi$ is a solution of the Euler-Lagrange equations $\Phi\spi{}\pb i_{\sigma\spi}\Omega_L=0$, then 
$$
\Phi\spi{}\pb d_{\sigma\spi}\L 
= d[\Phi\spi{}\pb i_{\sigma\spi}\Theta_L].
$$
The above is the base of Noether's theorem. 

\begin{definition}
We will say that the Lagrangian density $\L$ is \emph{invariant} under a $\pi$-vertical section $\sigma\in\sec{E}$ if $d_{\sigma\spi}\L=0$.
\end{definition}

It follows from this definition that $\L$ is invariant under $\sigma$ if and only if it is invariant under the flow of $\sigma\spi$, that is $\Phi_s\spi{}\pb\L=\L$. 

\begin{definition}
A $(r-1)$-form $\lambda$ is said to be a \emph{conserved current} if $\Phi\spi{}\pb\lambda$ is a closed form for any solution $\Phi$ of the Euler-Lagrange equations.
\end{definition}

With these definitions and taking into account the above arguments we have the following result. 
\begin{theorem}
Let $\sigma\in\sec{E}$ be a $\pi$-vertical section. If the Lagrangian density is invariant under $\sigma$ then $i_{\sigma\spi}\Theta_L$ is a conserved current.
\end{theorem}
\begin{proof}
If $d_{\sigma\spi}\L=0$, from the relation $\Phi\spi{}\pb d_{\sigma\spi}\L 
= d[\Phi\spi{}\pb i_{\sigma\spi}\Theta_L]$, we have that $d[\Phi\spi{}\pb i_{\sigma\spi}\Theta_L]$, for every solution $\Phi$ of the Euler-Lagrange equations.
\end{proof}

We insist however that the lack of an integral description does not allow to interpret this as a conservation law in the classical sense. 


\section{The Hamiltonian formalism}
\label{Hamiltonian}

For the Hamiltonian approach we consider the affine dual of $\Jpi$. This is the bundle whose fiber over $m\in M$ is the set of all affine maps from $\Jpi[m]$ to $\R$. As in the standard case, there is a more convenient representation of this bundle as a bundle of $r$-forms, which we will denote by the same symbol. Explicitly, we consider the bundle $\map{\prdJM}{\dJpi}{M}$ whose fibre at $m\in M$ is
$$
\dJpi[m]=\set{\lambda\in{(E^*_m)^{\wedge r}}}{i_{k_1}i_{k_2}\lambda=0 \text{ for all } k_1,k_2\in K_m}.
$$
An element in $\dJpi$ has the local expression $\lambda =\lambda_0\omega+\lambda^a_\alpha e^\alpha\wedge\omega_a$ and it is identified with the affine map $\lambda_0+\lambda^a_\alpha y^\alpha_a$. We thus have local coordinates $(x^i,u^A,\mu_0,\mu^a_\alpha)$ on $\dJpi$, where $\mu_0$ and $\mu^a_\alpha$ are the functions given by $\mu_0(\lambda)=\lambda_0$ and $\mu^a_\alpha(\lambda)=\lambda^a_\alpha$. We also consider the $E$-tangent bundle to $\dJpi$, whose fibre at the point $\lambda$ is
$$
\prol[E]{\dJpi}[\lambda]=
\set{(a,V)\in E_m\times T_\lambda\dJpi}{\rho_E(a)=T\prdJM(V)},
$$
and the projection $\prTdJE=\prol{\prdJM}$, which defines the bundle map $\pi_{10}^\dag=(\prTdJE,\prdJM)$.

On $\dJpi$ we have a canonically defined $r$-form $\Theta$, the \emph{multimomentum form}, given by $\Theta_\lambda=(\pi_{10}^\dag)\pb\lambda$, or explicitly,
$$
\Theta_\lambda(Z_1,Z_2,\ldots,Z_r)=\lambda(a_1,a_2,\ldots,a_r),
$$
for $Z_i=(\lambda,a_i,V_i)\in\prol[E]{\dJpi}[\lambda]$, $i=1,\ldots,r$. The canonical multimomentum form has the local expression
$$
\Theta=\mu_0\omega+\mu^a_\alpha\X^\alpha\wedge\omega_a.
$$
The differential of the canonical multimomentum form is (minus) the canonical \emph{multisymplectic form} on $\dJpi$, that is $\Omega=-d\Theta$. Its local expression is
$$
\Omega=
\X^\alpha\wedge\P^a_\alpha\wedge\omega_a
+\frac{1}{2}\mu^a_\gamma C^\gamma_{\alpha\beta}\X^\alpha\wedge\X^\beta\wedge\omega_a
-\P_0\wedge\omega
-\mu^a_\gamma(C^\gamma_{a\alpha}+C^b_{ab}\delta^\gamma_\alpha)\X^\alpha\wedge\omega.
$$

We also consider the vector bundle dual to $\Vpi$ represented as the bundle $\map{\prdVM}{\dVpi}{M}$ whose fibre at $m\in M$ is the set of all linear maps $\map{\xi}{K_m}{F_n^{\wedge(r-1)}}$, where $n=\prMN(m)$. An element of $\dVpi$ is of the form $\xi=\xi^a_\alpha e^\alpha\wedge\omega_a$, and therefore we have local coordinates $(x^i,u^A,\mu^a_\alpha)$ on $\dVpi$, where $\mu^a_\alpha(\xi)=\xi^a_\alpha$.

We have used the same symbol for the coordinates $\mu^a_\alpha$ in $\dJpi$ and in $\dVpi$ because the manifold $\dJpi$ is fibred over the manifold $\dVpi$. Indeed, an element $\lambda$ in $\dJpi[m]$ is an affine map and therefore has an associated linear map, which is represented as an element of $\dVpi$, and thus provides a surjective submersion $\map{\ell}{\dJpi}{\dVpi}$. As before, we consider the $E$-tangent bundle $\prol[E]{\dVpi}$ whose fibre at $\xi\in\dVpi$ is
$$
\prol[E]{\dVpi}[\xi] =
\set{(a,V)\in E_m\times T_\xi\dVpi}{\rho_E(a)=T\prdVM(V)}.
$$

In order to avoid confusions, the basis of local sections of $\prol[E]{\dJpi}$ will be denoted $(\X_a,\X_\alpha,\P^0,\P^\alpha_a)$ and similarly, the basis of local sections of $\prol[E]{\dVpi}$ will be denoted $(\X_a,\X_\alpha,\P^\alpha_a)$. That is, we will denote the vertical elements $\left(\xi,0,\pd{}{\mu^a_\alpha}\right)$ by $\P^\alpha_a$ instead of $\V^\alpha_a$ as before. Accordingly, the dual basis will be denoted $(\X^a,\X^\alpha,\P_0,\P_\alpha^a)$ and $(\X^a,\X^\alpha,\P_\alpha^a)$, respectively.

\subsection*{Liouville-Cartan forms}
By a Hamiltonian section we mean a section of the bundle $\map{\ell}{\dJpi}{\dVpi}$. With the help of a Hamiltonian section $h$ we can pullback the canonical multimomentum and multisymplectic forms to $\dVpi$,
$$
\Theta_h=(\prol{h})\pb\Theta
\qquand
\Omega_h=(\prol{h})\pb\Omega.
$$
Notice that $\Omega_h=-d\Theta_h$ since $\prol{h}$ is a morphism of Lie algebroids.
Locally, a Hamiltonian section $h$ is determined by a local function $H(x^i,u^A,\mu^a_\alpha)$ by means of $h(x^i,u^A,\mu^a_\alpha)=(x^i,u^A,\mu^a_\alpha,-H(x^i,u^A,\mu^a_\alpha))$, that is, the coordinate $\mu_0$ is determined by $\mu_0=-H(x^i,u^A,\mu^a_\alpha)$.
It follows that the forms $\Theta_h$ and $\Omega_h$ have the local expression
$$
\Theta_h=\mu^a_\alpha\X^\alpha\wedge\omega_a-H\omega
$$
and
$$
\Omega_h=
\X^\alpha\wedge\P^a_\alpha\wedge\omega_a
+\frac{1}{2}\mu^a_\gamma C^\gamma_{\alpha\beta}\X^\alpha\wedge\X^\beta\wedge\omega_a
+dH\wedge\omega
-\mu^a_\gamma(C^\gamma_{a\alpha}+C^b_{ab}\delta^\gamma_\alpha)\X^\alpha\wedge\omega.
$$

\subsection*{Hamilton equations}
A solution of the Hamiltonian system defined by a Hamiltonian section $h$ is a morphism $\Lambda$ from $\map{\tau^F_N}{F}{N}$ to $\map{\tau^E_{\dVpi}}{\prol[E]{\dVpi}}{\dVpi}$ such that
$$
\Lambda\pb\left(i_X\Omega_h\right)=0,
$$
for every section $X$ of $\prol[E]{\dVpi}$ vertical over $F$. The above equations are said to be the Hamiltonian field equations.

Such a map $\Lambda$ has the local expression $\Lambda=(\X_a + \Lambda^\alpha_a\X_\alpha + \Lambda^c_{\gamma a}\P^\gamma_a)\otimes\bar{e}^a$. By taking the section $X=\P^\alpha_a$ in the Hamilton equations we get
$$
i_X\Omega
=\pd{H}{\mu^a_\alpha}\omega-\X^\alpha\wedge\omega_a
$$
and then
$$
\Lambda\pb(i_X\Omega)=\left(\pd{H}{\mu^a_\alpha}-\Lambda^\alpha_a\right)\omega,
$$
from where we get
$$
\Lambda^\alpha_a=\pd{H}{\mu^a_\alpha}.
$$

Taking now the section $X=\X_\alpha$ in the Hamilton equations we get
$$
i_X\Omega
=\left(\rho^A_\alpha\pd{H}{u^A}
  -\mu^c_\gamma (C^\gamma_{c\alpha}+C^b_{ab}\delta^\gamma_\alpha)\right)\omega
+\left(\P^a_\alpha+\mu^a_\gamma C^\gamma_{\alpha\beta}\X^\beta\right)\wedge\omega_a
$$
and thus
$$
\Lambda\pb(i_X\Omega)=\left(
\rho^A_\alpha\pd{H}{u^A}-\mu^c_\gamma (C^\gamma_{c\alpha}+C^\gamma_{\beta\alpha}\Lambda^\beta_\alpha)
+\Lambda^c_{\alpha c}-\mu^c_\gamma C^b_{ab}
\right)\omega,
$$
from where we get the relation
$$
\Lambda^c_{\alpha c}=
-\rho^A_\alpha\pd{H}{u^A}+\mu^c_\gamma\left( C^\gamma_{c\alpha}+
C^\gamma_{\beta\alpha}\Lambda^\beta_c\right)-\mu^a_\beta C^b_{ab}.
$$

It follows that the Hamiltonian field equations are
$$
\Lambda^\alpha_a=\pd{H}{\mu^a_\alpha}
\qquand
\Lambda^c_{\alpha c}+\mu^b_\alpha C^c_{bc}=
-\rho^A_\alpha\pd{H}{u^A}+\mu^c_\gamma\left( C^\gamma_{c\alpha}+
C^\gamma_{\beta\alpha}\pd{H}{\mu^c_\beta}\right).
$$
Notice that they do not determine $\Lambda$, but they provide the value of $\Lambda^c_{\alpha c}$ the trace of $\Lambda^a_{\alpha b}$.

On the other hand, the map $\Lambda$ must be a morphism; in particular it is an admissible map $\rho^1\circ\ol{\Lambda}=T\ul{\Lambda}\circ\rho_F$, which in coordinates reads
$$
\rho^i_a\pd{u^A}{x^i}=\rho^A_a+\rho^A_\alpha\Lambda^\alpha_a
\qquand
\rho^i_a\pd{\mu^c_\gamma}{x^i}=\Lambda^c_{\gamma a}.
$$
The additional properties for being a morphism are
$$
 \left(\pd{H}{\mu^a_\alpha}\right)'_{|b}
-\left(\pd{H}{\mu^b_\alpha}\right)'_{|a}
+C^\alpha_{\beta\gamma}\pd{H}{\mu^b_\beta}\pd{H}{\mu^a_\gamma}
+C^\alpha_{b\gamma}\pd{H}{\mu^a_\gamma}
-C^\alpha_{a\gamma}\pd{H}{\mu^b_\gamma}
+C^\alpha_{ba}=0
$$
and
$$
\Lambda^a_{\alpha b|c}-\Lambda^a_{\alpha c|b}
=\Lambda^a_{\alpha d}C^d_{cb},
$$
but this one is a consequence of the admissibility conditions.
Therefore we get the system of partial differential equations
\begin{align*}
&\rho^i_a\pd{u^A}{x^i}=\rho^A_a+\rho^A_\alpha\pd{H}{\mu^a_\alpha}\\
&\left(\pd{H}{\mu^a_\alpha}\right)'_{|b}
-\left(\pd{H}{\mu^b_\alpha}\right)'_{|a}
+C^\alpha_{\beta\gamma}\pd{H}{\mu^b_\beta}\pd{H}{\mu^a_\gamma}
+C^\alpha_{b\gamma}\pd{H}{\mu^a_\gamma}
-C^\alpha_{a\gamma}\pd{H}{\mu^b_\gamma}
+C^\alpha_{ba}=0\\
&\rho^i_c\pd{\mu^c_\alpha}{x^i}+\mu^b_\alpha C^c_{bc}=
-\rho^A_\alpha\pd{H}{u^A}+\mu^c_\gamma\left( C^\gamma_{c\alpha}+
C^\gamma_{\beta\alpha}\pd{H}{\mu^c_\beta}\right).
\end{align*}
which will be called the Hamiltonian partial differential equations.

\subsection*{The Legendre transformation}
Let $\phi_0\in\Jpi$ and $m=\prJM(\phi_0)$. We define $\leg(\phi_0)\in\dJpi[m]$ as follows. For every $\phi\in\Jpi[m]$, we consider the map  $t\mapsto\varphi(t)=\L(\phi_0+t(\phi-\phi_0))$. The value of $\leg(\phi_0)$ over $\phi$ is the first order affine approximation of $\varphi$ at $t=0$,
$$
\leg(\phi_0)(\phi)=\varphi(0)+\varphi'(0).
$$
The map $\map{\leg}{\Jpi}{\dJpi}$ is said to be the Legendre map or the Legrendre transformation.

In coordinates, if $\phi_0=(e_a+y^\alpha_a e_\alpha)\otimes\bar{e}^a$ and $\phi=(e_a+z^\alpha_a e_\alpha)\otimes\bar{e}^a$, then $\varphi(t)=L(x^i,u^A,y^\alpha_a+t(z^\alpha_a-y^\alpha_a))$, so that
$$
\leg(\phi_0)(\phi)=\varphi(0)+\varphi'(0)
=L(x^i,u^A,y^\alpha_a)
 +\pd{L}{y^\alpha_a}(x^i,u^A,y^\alpha_a)(z^\alpha_a-y^\alpha_a),
$$
which under the canonical identification $(p+p^a_\alpha z^\alpha_a)\equiv p\omega+p^a_\alpha e^\alpha\wedge\omega_a$ corresponds to 
$$
\leg(\phi_0)(\phi)
=
\pd{L}{y^\alpha_a}(e^\alpha-y^\alpha_be^b)\wedge\omega_a+
\left(L-\pd{L}{y^\alpha_a}y^\alpha_a\right)\omega,
$$
It follows from this expression that $\leg$ is smooth and that $(\prol{\leg})\pb\Theta=\Theta_\L$. Therefore $(\prol{\leg})\pb\Omega=\Omega_\L$

Finally, the reduced Legendre transformation is the map $\map{\legd}{\Jpi}{\dVpi}$ given by projection $\legd=\ell\circ\leg$. In coordinates, if $\psi\in\Vpi$ has coordinates $\psi=(x^i,u^A,v^\alpha_a)$ then we have 
$$
\legd(\phi_0)(\psi)=\pd{L}{y^\alpha_a}(x^i,u^A,y^\alpha_a)\,v^\alpha_a,
$$
which under the canonical identification  $p^a_\alpha z^\alpha_a\equiv p^a_\alpha e^\alpha\wedge\omega_a$ corresponds to
$$
\legd(\phi_0)(\psi)
=
\pd{L}{y^\alpha_a}e^\alpha\wedge\omega_a.
$$

In local coordinates the tangent prolongation of the Legendre transformation is given by
\begin{align*}
&(\prol{\leg})\pb x^i=x^i
&&(\prol{\leg})\pb \X^a=\X^a\\
&(\prol{\leg})\pb u^A=u^A
&&(\prol{\leg})\pb \X^\alpha=\X^\alpha\\
&(\prol{\leg})\pb\mu_0=L-\pd{L}{y^\alpha}y^\alpha_a\equiv E_L
&&(\prol{\leg})\pb \P_0=dE_L\\
&(\prol{\leg})\pb\mu_\alpha^a=\pd{L}{y^\alpha}
&&(\prol{\leg})\pb \P^a_\alpha=d\left(\pd{L}{y^\alpha}\right)
\end{align*}
and 
\begin{align*}
&(\prol{\legd})\pb x^i=x^i
&&(\prol{\legd})\pb \X^a=\X^a\\
&(\prol{\legd})\pb u^A=u^A
&&(\prol{\legd})\pb \X^\alpha=\X^\alpha\\
&(\prol{\legd})\pb\mu_\alpha^a=\pd{L}{y^\alpha}
&&(\prol{\legd})\pb \P^a_\alpha=d\left(\pd{L}{y^\alpha}\right)
\end{align*}

\subsection*{Equivalence}
If the Legendre transformation $\legd$ is a (global) diffeomorphism we will say that the Lagrangian is hyperregular. In this case both formalisms, the Lagrangian and the Hamiltonian, are equivalent. 

\begin{proposition}
The following conditions are equivalent:
\begin{enumerate}
\item $\legd$ is a local diffeomorphism.
\item For every $\phi\in\Jpi$ the linear map $\prol{\legd}[\phi]$ is invertible.
\item The Lagrangian $L$ is regular.
\end{enumerate}
\end{proposition}
\begin{proof}{}
 [$(1)\Leftrightarrow (2)$] If $\legd$ is a local diffeomorphism, then $T\legd$ is invertible. Thus, for every $\phi\in\Jpi$ we have that  $(\legd(\phi),a,V)\mapsto (\phi,a,(T_\phi\legd)^{-1}(V))$ is the inverse of $\prol{\legd}[\phi]$. Conversely, if $\prol{\legd}[\phi]$ is invertible for every $\phi\in\Jpi$, and $V\in T_\phi\Jpi$ is such that $T_\Phi\legd(V)=0$, then $V$ is vertical (because $\legd$ is a map over the identity) and thus 
$$
\prol{\legd}[\phi](\phi,0,V)=(\legd(\phi),0,T_\phi\legd(V))=(\legd(\phi),0,0).
$$
Since $\prol{\legd}[\phi]$ in invertible, we get that $V=0$. 

[$(2)\Leftrightarrow (3)$] From the local expression of the prolongation of $\legd$ we get that $\prol{\legd}[\phi]$ is invertible if and only if the matrix $\pd{^2L}{y^\alpha_a\partial y^\beta_b}$ is invertible, which is equivalent to the regularity of the Lagrangian.
\end{proof}

\begin{theorem}
Let $L$ be a hyperregular Lagrangian. If $\Phi$ is a solution of the Euler-Lagrange equations then $\Lambda=\prol{\legd}\circ\Phi\spi$ is a solution of the Hamiltonian field equations. Conversely, if $\Lambda$ is a solution of the Hamiltonian field equations then there exists one and only one solution $\Phi$ of the Euler-Lagrange equations such that $\Lambda=\prol{\legd}\circ\Phi\spi$.
\end{theorem}
\begin{proof}
Let $\Phi\in\sec{\pi}$ be a solution of the Euler-Lagrange equations and define $\Lambda=\prol{\legd}\circ\Phi\spi$. Then $\Lambda$ is a morphism, because so is $\Phi\spi$ and $\prol{\legd}$. Moreover, it is a solution of the Hamiltonian field equations. Indeed, if $X$ and $Y$ are such that $\prol{\legd}\circ Y=X\circ\legd$, then 
$$
\Lambda\pb i_X\Omega
=\Phi\spi{}\pb(\prol{\legd})\pb i_X\Omega
=\Phi\spi{}\pb i_Y(\prol{\legd})\pb\Omega
=\Phi\spi i_Y\Omega_L.
$$
If $\Phi$ is a solution of the Euler-Lagrange equations then the right hand side vanishes for every $Y$ vertical, and then the left hand side vanishes for every $X$ vertical. (Notice that if $Y$ is vertical and arbitrary so is $X$.)

Conversely, let $\Lambda$ be a solution of the Hamiltonian field equations, and define $\Psi=(\prol{\legd})^{-1}\circ\Lambda$. Then, for every vertical $Y$ we have 
$$
\Psi\pb i_Y\Omega_L
=\Psi\pb i_Y(\prol{\legd})\pb\Omega
=\Psi\pb (\prol{\legd})\pb i_X\Omega
=(\prol{\legd}\circ\Psi)\pb i_X\Omega
=\Lambda\pb i_X\Omega.
$$
so that $\Psi$ satisfies $\Psi\pb i_X\Omega_L=0$, for every vertical $X$. Since the Lagrangian is regular we have that $\Psi$ is semiholonomic. But since $\Lambda$ is a morphism, so is $\Psi$, and hence $\Psi$ is jet prolongation $\Psi=\Phi\spi$ with $\Phi$ a morphism. Therefore $\Phi\in\sec{\pi}$ is a solution of the Euler-Lagrange field equations.
\end{proof}

Unfortunately, the most interesting examples of Lagrangian field theories are singular, so that the above result does not apply. In such cases we can proceed as in~\cite{LeMSMa} or alternatively as it is explained in the next subsection.

\subsection*{The unified formalism}

Probably the best alternative in the analysis of the Hamiltonian equations is the so-called unified formalism~\cite{Barna-U, LeMaMD}.

Consider the fibred product ${\Jpi\times_M\dVpi}\to{M}$, and the map $\map{\lambda_L}{\Jpi\times_M\dVpi}{\dJpi}$ defined by 
$$
\lambda_L(\phi,\psi)(\phi')=L(\phi)+\psi(\phi'-\phi).
$$
In local coordinates, if $\phi=(x^i,u^A,y^\alpha_a)$, $\psi=(x^i,u^A,\mu^a_\alpha)$ and $\phi'=(x^i,u^A,z^\alpha_a)$, then 
$$
\lambda_L(\phi,\psi)(\phi')=L(x^i,u^A,y^\alpha_a)+\mu^a_\alpha(z^\alpha_a-y^\alpha_a).
$$
Thus the coordinates of $\lambda_L(\phi,\psi)$ are $(x^i,u^A,L(x^i,u^A,y^\alpha_a)-\mu^a_\alpha y^\alpha_a,\mu^a_\alpha)$, that is the value of the coordinate $\mu_0$ of $\lambda_L(\phi,\psi)$ is $\mu_0=L(x^i,u^A,y^\alpha_a)-\mu^a_\alpha y^\alpha_a$.

We consider the $E$-tangent bundle $\map{\tau^E_{\Jpi\times_M\dVpi}}{\prol[E]{(\Jpi\times_M\dVpi)}}{\Jpi\times_M\dVpi}$ and the projections $\upsilon_{10}$  and $\upsilon_1$ from this bundle to $\map{\tau^E_M}{E}{M}$ and $\map{\tau^F_N}{F}{N}$, respectively.

By pulling back the canonical multisymplectic form $\Omega$ on $\dJpi$ to $\Jpi\times_M\dVpi$ by $\prol{\lambda_L}$ we get a multisymplectic form $\Omega^L$ on $\Jpi\times_M\dVpi$,
$$
\Omega^L=(\prol{\lambda_L})\pb\Omega,
$$
in terms of which we set the field equations. Notice that if we define $\Theta^L=(\prol{\lambda_L})\pb\Theta$ then $\Omega^L=-d\Theta^L$. In coordinates 
$$
\Theta^L=(L-\mu_\alpha^a y_a^\alpha)\omega +\mu^a_\alpha\X^\alpha\wedge\omega_a
=L\omega +\mu^a_\alpha\theta^\alpha\wedge\omega_a.
$$

\begin{definition}
By a  solution of unified field equations we mean a morphism $\Upsilon\in\sec{\upsilon_1}$ such that $\Upsilon\pb(i_X\Omega^L)=0$ for every $\upsilon_1$-vertical section $X$ of $\prol[E]{(\Jpi\times_M\dVpi)}$.
\end{definition}

We will see that this equations reproduce the field equations and the semiholonomy condition. 

The section $\Upsilon$ is of the form $\Upsilon=(\Psi,\Lambda)$, where $\Psi$ is a section of $\Jpi$ and $\Lambda$ is a section of $\dVpi$, both projecting to the same section of $\pi$. It is easy to see that $\Upsilon$ is a morphism of Lie algebroids if and only if $\Psi$ and $\Lambda$ are morphisms of Lie algebroids.

We consider the forms $\Delta^L=\Theta^L-(\prol{\pr_1})\pb\Theta_L$ and $\Xi^L=\Omega^L-(\prol{\pr_1})\pb\Omega_L$, so that $\Xi^L=-d\Delta^L$. In coordinates we have that 
$$
\Delta^L=\left(\mu^a_\alpha-\pd{L}{y^\alpha_a}\right)\theta^\alpha\wedge\omega_a.
$$
\begin{proposition}
Let $\Upsilon\in\sec{\upsilon_1}$ be a morphism, and $\Psi,\Lambda$ its components.
\begin{enumerate}
\item $\Upsilon\pb i_Z\Omega^L=0$ for every $\prol{\pr_1}$-vertical section $Z$ if and only if $\Psi$ is semiholonomic (and hence holonomic).
\item $\Upsilon\pb i_Z\Omega^L=0$ for every $\prol{\upsilon}$-vertical section $Z$ if and only if $\Psi$ is semiholonomic and $\textrm{Im}\ul{\Upsilon}\subset\mathrm{Graph}(\legd)$.
\item $\Upsilon\pb i_Z\Omega^L=0$ for every $\prol{\upsilon}$-vertical section $Z$ if and only if $\Upsilon\pb i_W\Omega^L=0$ for every section $W$.
\end{enumerate}
\end{proposition}
\begin{proof}
Indeed, for $Z=\P^\alpha_a$ we have that 
$$
i_{\P^\alpha_a}\Xi^L=-\theta^\alpha\wedge\omega_a,
$$
so that $\Upsilon\pb i_{\P^\alpha_a}\Xi^L = \Psi\pb\theta^\alpha\wedge\omega$ vanishes if and only if $\Psi$ is semiholonomic. Since $\Psi$ is a morphism we have that $\Psi$ is holonomic. This proves the first condition.

For the second we have 
$$
i_{\V_\alpha^a}\Xi^L=\pd{^2L}{y^\alpha_a\partial y^\beta_b}\theta^\beta\wedge\omega_b+\left(\mu^a_\alpha-\pd{L}{y^\alpha_a}\right)\omega.
$$
Therefore $\Upsilon\pb i_Z\Omega^L=0$ for $Z=\P^\alpha_a$ and $Z=\V_\alpha^a$ if and only if $\Psi$ is semiholonomic and $\Upsilon\pb\left(\mu^a_\alpha-\pd{L}{y^\alpha_a}\right)=0$, that is, $\ul{\Upsilon}$ takes values on the graph of the Legendre transformation $\legd$.

Finally, for the third we take $Z=\X_\alpha$ and thus
$$
i_{\X_\alpha}\Xi^L
=-d_{\X_\alpha}\Delta^L+d(i_{\X_\alpha}\Delta_L).
$$
The pullback by $\Upsilon$  of the second term vanishes because $\Upsilon\pb$ commutes with $d$ and 
$$
\Upsilon\pb i_{\X_\alpha}\Delta_L=\Upsilon\pb\left(\mu^a_\alpha-\pd{L}{y^\alpha_a}\right)\omega_a=0.
$$
Thus we have $\Upsilon\pb i_{\X_\alpha}\Xi^L
=-\Upsilon\pb d_{\X_\alpha}\Delta^L$ which vanishes as a consequence of the first two properties
\begin{align*}
\Upsilon\pb d_{\X_\alpha}\Delta^L
&=\Upsilon\pb\left[d_{\X_\alpha}\left(\mu^b_\beta-\pd{L}{y^\beta_b}\right)\right]
\left[\Upsilon\pb\theta^\beta\wedge\omega_a\right]+\\
&\qquad\qquad{}+
\Upsilon\pb\left(\mu^b_\beta-\pd{L}{y^\beta_b}\right)
\Upsilon\pb\left[d_{\X_\alpha}(\theta^\beta\wedge\omega_b)\right]\\
&=0+0=0.
\end{align*}
This proves the third.
\end{proof}

We notice that if $\Psi$ is semiholonomic and $\Upsilon=(\Psi,\Lambda)$ is a morphism of Lie algebroids (admissible is sufficient) such that the base map is in the graph of the Legendre transformation,  $\textrm{Im}\ul{\Upsilon}\subset\mathrm{Graph}(\legd)$, then $\Upsilon$ itself is in the graph of $\prol{\legd}$, that is,  $\Lambda=\prol{\legd}\circ\Psi$. Indeed, if we denote by $\zeta$ the bundle map $\zeta=\Upsilon-(\Psi,\prol{\legd}\circ\Psi)$, then we have that $\zeta\pb x^i=0$, $\zeta\pb u^A=0$, $\zeta\pb y^\alpha_a=0$ and 
$$
\zeta\pb \mu^a_\alpha=\left(\mu^a_\alpha-\pd{L}{y^\alpha_a}\right) \circ\ul{\Upsilon}
$$
which vanishes if $\textrm{Im}\ul{\Upsilon}\subset\mathrm{Graph}(\legd)$. Moreover we have that $\zeta\pb\X^a=0$, $\zeta\pb\X^\alpha=0$ (because $\Psi$ is semiholonomic), $\zeta\pb\V^\alpha_a=0$, and
$$
\zeta\pb\P^a_\alpha
=\zeta\pb d\mu^a_\alpha
=d\zeta\pb\mu^a_\alpha
=d\left[\left(\mu^a_\alpha-\pd{L}{y^\alpha_a}\right) \circ\ul{\Upsilon}\right]
$$
which vanishes because $\textrm{Im}\ul{\Upsilon}\subset\mathrm{Graph}(\legd)$.
\begin{proposition}
Let $\Upsilon=(\Psi,\Lambda)$ be a solution of the unified field equations. Then $\Psi=\Phi\spi$ with $\Phi$ a solution of the Euler-Lagrange equations and $\Lambda=\prol{\legd}\circ\Phi\spi$. 

Conversely, if $\Phi$ is solution of the Euler-Lagrange equations then $\Upsilon=(\Phi\spi,\prol{\legd}\circ\Phi\spi)$ is a solution of the unified field equations.
\end{proposition}
\begin{proof}
Let $\Upsilon$ be a solution of the unified field equations. Then if $Z=(X,Y)$ is a section of $\prol[E]{\Jpi\times_M\dVpi}$ we have
$$
\Psi\pb i_X\Omega_L=\Upsilon\pb i_Z\Xi^L.
$$
Thus for a $\prol{\pr_1}$-vertical $Z$ we have that $X=0$, and thus $\Upsilon\pb i_Z\Xi^L=0$, so that $\Psi$ is semiholonomic, and hence holonomic $\Psi=\Phi\spi$. On the other hand if  $Z$ is $\prol{\pr_2}$-vertical, then $X$ is vertical and by proposition~\ref{semiholonomic.and.vertical} we get that $\Psi\pb i_X\Omega_L=0$. Thus $\Upsilon\pb i_Z\Xi^L=0$ and it follows that $\Upsilon$ is in the graph of the Legendre transformation, i.e. it is of the form $\Upsilon=(\Psi,\prol{\legd}\circ\Psi)$. Finally, since $\Upsilon\pb i_Z\Xi^L$ vanishes for verticals $Z$, it vanishes for every $Z$, and hence $\Psi\pb i_X\Omega_L=0$. Thus $\Psi=\Phi\spi$ satisfies the De Donder equations, i.e. $\Phi$ is a solution of the Euler-Lagrange equations.

Conversely, let $\Phi$ be a solution of the Euler-Lagrange equations and set $\Upsilon=(\Phi\spi,\prol{\legd}\circ\Phi\spi)$. Since $\Phi$ is a morphism so is $\Phi\spi$, hence so is $\prol{\legd}\circ\Phi\spi$ and hence so is $\Upsilon$. Thus
$$
\Upsilon\pb i_Z\Omega^L=\Psi\pb i_X\Omega_L+\Upsilon\pb i_Z\Xi_L=\Upsilon\pb i_Z\Xi_L,
$$
which vanishes because $\Phi\spi$ is semiholonomic and (by definition) the image of $\Upsilon$ is contained in the graph of the Legendre transformation.
\end{proof}


\section{Examples}
\label{examples}

\subsection*{Standard case} In the standard case, we consider a bundle $\map{\ul{\pi}}{M}{N}$, the standard Lie algebroids $F=TN$ and $E=TM$ and the tangent map $\map{\ol{\pi}=T\ul{\pi}}{TM}{TN}$. Then we have that $\Jpi=J^1\ul{\pi}$. Thus we recover the standard first order Field Theory defined on sections on the bundle $M$. If we take coordinate basis, we recover the local expression of the Euler-Lagrange field equations and the Hamiltonian field equations. Moreover, if we take general non-coordinate basis of vectorfields, our equations provides an expression of the standard Euler-Lagrange and Hamiltonian field equations written in pseudo-coordinates.

In particular, one can take an Ehresmann connection on the bundle $\map{\ul{\pi}}{M}{N}$ and use an adapted local basis 
$$
\bar{e}_i=\pd{}{x^i}
\qquand
\left\{\begin{aligned}
&e_i=\pd{}{x^i}+\Gamma^A_i\pd{}{u^A}\\
&e_A=\pd{}{u^A}.
\end{aligned}\right.
$$
In this case, the greek indices and the latin capital indices coincide. Then we have the brackets 
$$
[e_i,e_j]=-R^A_{ij}e_A,
\qquad
[e_i,e_B]=\Gamma^A_{iB}e_A
\qquand
[e_A,e_B]=0,
$$
where we have written $\Gamma^B_{iA} = \partial{\Gamma^B_i}/\partial{u^A}$ and where $R^A_{ij}$ is the curvature tensor of the nonlinear connection we have chosen. The components of the anchor are $\rho^i_j=\delta^i_j$, $\rho^A_i=\Gamma^A_i$ and $\rho^A_B=\delta^A_B$ so that the Euler-Lagrange equations are 
\begin{align*}
&\pd{u^A}{x^i}=\Gamma^A_i+y^A_i\\
&\pd{y^A_{i}}{x^j}-\pd{y^A_{j}}{x^i}
+\Gamma^A_{jB}y^B_i-\Gamma^A_{iB}y^B_j
=R^A_{ij}\\
&\frac{d}{dx^i}\left(\pd{L}{y^A_i}\right)
-\Gamma^B_{iA}\pd{L}{y^B_i}
=\pd{L}{u^A}.
\end{align*} 
On the Hamiltonian side we have the Hamiltonian field equations
\begin{align*}
&\pd{u^A}{x^i}=\Gamma^A_i+\pd{H}{\mu^i_A}\\
&\frac{d}{dx^j}\left(\pd{H}{\mu^i_A}\right)
-\frac{d}{dx^i}\left(\pd{H}{\mu^j_A}\right)
+\Gamma^A_{jB}\pd{H}{\mu^i_B}
-\Gamma^A_{iB}\pd{H}{\mu^j_B}
=R^A_{ij}\\
&\pd{\mu^i_A}{x^i}-\mu^i_B\Gamma^B_{iA}=-\pd{H}{u^A}.
\end{align*}
See~\cite{CaCrIb,Barna-L,Barna-H} for alternative derivations of this equations.

\subsection*{Time-dependent Mechanics} In~\cite{SaMeMa1,SaMeMa2} we developed a theory of Lagrangian and Hamiltonian mechanics for time dependent systems defined on Lie algebroids, where the base manifold is fibred over the real line $\R$. Later, in~\cite{MaMeSa}, we generalize such results to the case of a general manifold (not necessarily fibred over $\R$) which was based on the notion of Lie algebroid structure over an affine bundle. Since time-dependent mechanics is but a 1-dimensional Field Theory, our results must be related to that. 

The case where the manifold $M$ is fibred over $\R$ fits in the scheme developed in this paper as follows. Consider a Lie algebroid $\map{\tau^E_M}{E}{M}$ and the standard Lie algebroid $\map{\tau_\R}{T\R}{\R}$. We consider the Lie subalgebroid $K=\ker(\pi)$ and define 
$$
A=\set{a\in E}{\ol{\pi}(a)=\pd{}{t}}.
$$
Then $A$ is an affine subbundle modeled on $K$ and the `bidual' of $A$ is $(A^\dag)^*=E$. Moreover, the Lie algebroid structure on $E$ defines by restriction a Lie algebroid structure on the affine bundle $A$.

Conversely, let $A$ be an affine bundle with a Lie algebroid structure as in~\cite{MaMeSa}. Then the vector bundle $E\equiv(A^\dag)^*$ has an induced Lie algebroid structure. If $\tilde{\rho}$ is the anchor of this bundle then the map $\ol{\pi}$ defined by $\ol{\pi}(z)=T\pi(\tilde{\rho}(z))$ is a morphism. Moreover we have that $A=\set{a\in E}{\ol{\pi}(a)=\pd{}{t}}$ as above.

We have a canonical identification of $A$ with $\Jpi$. Indeed, let $\map{T}{\Jpi}{A}$ be the map $T(\phi)=\Phi(\pd{}{t})$. Then $T$ is well defined, i.e. $T(\phi)\in A$ because
$$
\ol{\pi}(a)=\ol{\pi}(\phi(\pd{}{t}))=(\ol{\pi}\circ\phi)(\pd{}{t}) =\pd{}{t}.
$$
The inverse of the map $T$ is clearly $T^{-1}(a)=a\,dt$, where by $\phi=a\,dt$ we mean the map $\phi(\tau\pd{}{t})=\tau a$. One can easily follow the correspondences between both theories. For instance the vertical endomorphism $S_{dt}$ associated to the volume form $dt$ is the vertical endomorphism as defined in~\cite{MaMeSa}.

The morphism condition is just the admissibility condition so that the Euler-Lagrange equations are
$$
\begin{aligned}
&\frac{du^A}{dt}=\rho^A_0+\rho^A_\alpha y^\alpha\\
&\frac{d}{dt}\left(\pd{L}{y^\alpha}\right)
=\pd{L}{y^\gamma}(C^\gamma_{0\alpha}+C^\gamma_{\beta\alpha}y^\beta)
+\pd{L}{u^A}\rho^A_\alpha,
\end{aligned}
$$
where we have written $x^0\equiv t$ and $y^\alpha_0\equiv y^\alpha$. This expression is in agreement with~\cite{MaMeSa}. For the Hamilton equations we have
\begin{align*}
&\frac{du^A}{dt}=\rho^A_0+\rho^A_\alpha\pd{H}{\mu_\alpha}\\
&\frac{d\mu_\alpha}{dt}
=\mu_\gamma(C^\gamma_{0\alpha}+C^\gamma_{\beta\alpha}\pd{H}{\mu_\beta})
-\rho^A_\alpha\pd{H}{u^A}.
\end{align*}
where we have written $\mu^0_\alpha\equiv \mu_\alpha$.

\subsection*{The autonomous case}

A very common case in applications is that of systems depending on morphisms from a Lie algebroid to another one, which are not sections of any bundle. In this case we construct the bundle $E$ as the product of both algebroids. We assume that we have two Lie algebroids $\map{\tau^F_N}{F}{N}$ and $\map{\tau^G_Q}{G}{Q}$ over different bases and we set $M=N\times Q$ and $E=F\times G$, where the projections are both the projection over the first factor $\ul{\pi}(n,q)=n$ and $\ol{\pi}(a,k)=a$. The anchor is the sum of the anchors 
and the bracket is determined by the brackets of sections of $F$ and $G$ (a section of $F$ commutes with a section of $G$). We therefore have that 
$$
\rho^\alpha_a=0,\qquad C^\alpha_{ab}=0 \qquand C^\alpha_{a\beta}=0.
$$

In this case a jet at a point $(n,q)$ is of the form $\phi(a)=(a,\zeta(a))$, so that it is equivalent to a map $\map{\zeta}{F_n}{G_q}$. Therefore, in what follows we identify $\Jpi$ with the set of linear maps from a fibre of $F$ to a fibre of $G$.
This is further justified by the fact that a map $\map{\Phi}{F}{G}$ is a morphism of Lie algebroids if and only if the section $\map{(\id,\Phi)}{F}{F\times G}$ of $\pi$ is a morphism of Lie algebroids.

The affine functions $Z^\gamma_{a\alpha}$ reduce to  $Z^\gamma_{a\alpha}=C^\gamma_{\beta\alpha}y^\beta_a$ and thus the Euler-Lagrange equations are
$$
 \left(\pd{L}{y^\alpha_a}\right)'_{|a}
+C^b_{ba}\left(\pd{L}{y^\alpha_a}\right)
=\pd{L}{y^\gamma_a}C^\gamma_{\beta\alpha}y^\beta_a
+\pd{L}{u^A}\rho^A_\alpha.
$$

In the more particular and common case where $F=TN$ we can take a coordinate basis, so that we also have $C^c_{ab}=0$. Therefore the Euler-Lagrange partial differential equations are
$$
\begin{aligned}
&\pd{u^A}{x^a}=\rho^A_\alpha y^\alpha_a\\
&\frac{d}{dx^a}\left(\pd{L}{y^\alpha_a}\right)
=\pd{L}{y^\gamma_a}C^\gamma_{\beta\alpha}y^\beta_a
+\pd{L}{u^A}\rho^A_\alpha,\\
&\pd{y^\alpha_a}{x^b}-\pd{y^\alpha_b}{x^a}+C^\alpha_{\beta\gamma}y^\beta_by^\gamma_a=0,
\end{aligned}
$$

\subsubsection*{Autonomous Classical Mechanics} When moreover $F=T\R\to \R$ then we recover Weinstein's equations for a Lagrangian system on a Lie algebroid
$$
\begin{aligned}
&\frac{du^A}{dt}=\rho^A_\alpha y^\alpha\\
&\frac{d}{dt}\left(\pd{L}{y^\alpha}\right)
=\pd{L}{y^\gamma}C^\gamma_{\beta\alpha}y^\beta
+\pd{L}{u^A}\rho^A_\alpha,
\end{aligned}
$$
where, as before, we have written $x^0\equiv t$ and $y^\alpha_0\equiv y^\alpha$.  For the Hamilton equations we have 
\begin{align*}
&\frac{du^A}{dt}=\rho^A_\alpha\pd{H}{\mu_\alpha}\\
&\frac{d\mu_\alpha}{dt}
=\mu_\gamma C^\gamma_{\beta\alpha}\pd{H}{\mu_\beta}
-\rho^A_\alpha\pd{H}{u^A}.
\end{align*}
where we have written $\mu_\alpha^0\equiv \mu^0_\alpha$.

\subsection*{Poisson Sigma model} 

As an example of autonomous theory, we consider a 2-dimensional manifold $N$ and it tangent bundle $F=TN$. On the other hand, consider a Poisson manifold $(Q,\Lambda)$. Then the cotangent bundle $G=T^*Q$ has a Lie algebroid structure, where the anchor is $\rho(\sigma)=\Lambda(\sigma,\,\cdot\,)$ and the bracket is $[\sigma,\eta]=d^{TQ}_{\rho(\sigma)}\eta-d^{TQ}_{\rho(\eta)}\sigma-d^{TQ}\Lambda(\sigma,\eta)$, where $d^{TQ}$ is the ordinary exterior differential on~$Q$. 

The Lagrangian density for the Poisson Sigma model is $\L(\phi)=-\frac{1}{2}\phi\pb\Lambda$. In coordinates $(x^1,x^2)$ on $N$ and $(u^A)$ in $Q$ we have that $\Lambda=\frac{1}{2}\Lambda^{JK}\pd{}{u^J}\wedge\pd{}{u^K}$. 
A jet at the point $(n,q)$ is a map $\map{\phi}{T_nN}{T^*_qQ}$, locally given by $\phi=y_{Ki}du^K\otimes dx^i$. Thus we have local coordinates $(x^i,u^K,y_{Ki})$ on $\Jpi$. The local expression of the Lagrangian density is
$$
\L=-\frac{1}{2}\Lambda^{JK}A_J\wedge A_K=-\frac{1}{2}\Lambda^{JK}y_{J1}y_{K2}\,dx^1\wedge dx^2.
$$
where we have written $A_K=\Phi\pb(\partial/\partial u^K)=y_{Ki}dx^i$.

A long but straightforward calculation shows that for the Euler-Lagrange equation
$$
\frac{d}{dx^a}\left(\pd{L}{y^\alpha_a}\right)
=\pd{L}{y^\gamma_a}C^\gamma_{\beta\alpha}y^\beta_a
+\pd{L}{u^A}\rho^A_\alpha
$$
the right hand side vanishes while the left hand side reduces to 
$$
\frac{1}{2}\Lambda^{LJ}\left(
y_{L2|1}-y_{L1|2}+\pd{\Lambda^{MK}}{u^L}y_{M1}y_{K2}
\right)=0.
$$
In view of the morphism condition, we see that this equation vanishes. Thus the field equations are just
\begin{align*}
&\pd{u^J}{x^a}+\Lambda^{JK}y_{Ka}=0\\
&\pd{y_{Ja}}{x^b}-\pd{y_{Jb}}{x^a}+\pd{\Lambda^{KL}}{u^J} y_{Kb}y_{La}=0,
\end{align*}
or in other words
\begin{align*}
&d\phi^J+\Lambda^{JK}A_K=0\\
&d{A_J}+\frac{1}{2}\Lambda^{KL}_{,J} A_{K}\wedge A_{L}=0.
\end{align*}

The conventional Lagrangian density for the Poisson Sigma model~\cite{Strobl} is
$
\L'=\tr(\ol{\Phi}\wedge T\ul{\Phi})+\frac{1}{2}\Phi\pb\Lambda,
$
which in coordinates reads
$
\L'=A_J\wedge d\phi^J+\frac{1}{2}\Lambda^{JK}A_J\wedge A_K.
$
The difference between $\L'$ and $\L$ is a multiple of the admissibility condition $d\phi^J+\Lambda^{JK}A_K$;
$$
\L'-\L=
A_J\wedge d\phi^J+\frac{1}{2}\Lambda^{JK}A_J\wedge A_K
+\frac{1}{2}\Lambda^{JK}A_J\wedge A_K
=A_J\wedge(d\phi^J+\Lambda^{JK}A_K).
$$
Therefore both Lagrangians coincide on admissible maps, and hence on morphisms, so that the actions defined by them are equal.

\begin{remark} 
In the conventional analysis of the Poisson Sigma model~\cite{Strobl} the morphism condition is not imposed a priori and it is a result of the field equations. The 1-form $A_J$ acts as a Lagrange multiplier and something special occurs for this model since the multiplier is known in advance, not as consequence of the field equations.
\end{remark}

\begin{remark}
In more generality, one can consider a presymplectic Lie algebroid, that is, a Lie algebroid with a 2-cocycle $\Omega$, and the Lagrangian density $L=-\frac{1}{2}\Phi\pb\Omega$. The Euler-Lagrange equations vanishes as a consequence of the morphism condition and the closure of $\Omega$ so that we again get a topological theory. In this way one can generalize the theory for Poisson structures to a theory for Dirac structures.
\end{remark}

\subsection*{Holomorphic maps} 

Given two complex manifolds $(N,J_N)$ and $(M,J_M)$ a map $\map{\varphi}{N}{M}$ is a holomorphic map if its tangent map commutes the complex structures, that is, $T\varphi\circ J_N = J_M\circ T\varphi$. In many applications one has to consider a Lagrangian depending on holomorphic sections of a bundle $\map{\nu}{M}{N}$, the projection $\nu$ being holomorphic $T\nu\circ J_M=J_N\circ T\nu$. Our theory can include also such systems as follows. 

Since a complex structure $J$ on a manifold $Q$ satisfies $d_J^2=0$, we have that $TQ$ can be endowed with a Lie algebroid structure in which the exterior differential is $d_J=[i_J,d^{TQ}]$. Therefore a first idea is to use such Lie algebroid structures on both complex manifolds. Nevertheless, we have that a bundle map $\Phi=(\ol{\Phi},\ul{\Phi})$ is a morphism if and only if it is admissible, i.e. $J_M\circ\ol{\Phi}=T\ul{\Phi}\circ J_N$, which determines the map $\ol{\Phi}$ by $\ol{\Phi}=-J_M\circ T\ul{\Phi}\circ J_N$. Therefore, we do not get a holomorphic maps except if we impose the additional condition $\ol{\Phi}=T\ul{\Phi}$.

Since tangent maps are the only morphisms between tangent bundles (with the standard Lie algebroid structure), we can solve the problem by considering together both Lie algebroid structures on a complex manifold $Q$,
\begin{itemize}
\item the standard one, where the exterior differential is the standard exterior differential $d^{TQ}$,and
\item the Lie algebroid structure provided by the complex structure, where the exterior differential is $d_J=[i_J,d^{TQ}]$.
\end{itemize}
Notice that both structures are compatible in the sense that $[d_J,d^{TQ}]=0$. 

This is done by introducing an additional parameter $\lambda\in\R$ as follows. For a complex manifold $(Q,J)$ we set $G=\R\times TQ\to\R\times Q$, the vector bundle with projection $\tau^G_{\R\times Q}(\lambda,v)=(\lambda,\tau_Q(v))$. On $G$ we consider the Lie algebroid structure given by the exterior differential $d=d^{TQ}+\lambda d_{J}$, with $d\lambda=0$. Then it is clear that $d^2=0$ so that we have endowed $G$ with a Lie algebroid structure. The anchor is $\rho_G=(0,\id_{TQ}+\lambda J)$, that is 
$$
\rho_G(\lambda,v)=(0,v+\lambda Jv)\in T\R\times TQ,
$$

We consider the above construction for the complex manifold $N$ and we get the Lie algebroid $F=\R\times N$, and for the complex manifold $M$ and we get the Lie algebroid $E=\R\times M$. The projection $\pi=(\ol{\pi},\ul{\pi})$ given by $\ol{\pi}(\lambda,v)=(\lambda,T\nu(v))$ and $\ul{\pi}(\lambda,m)=(\lambda,\nu(m))$ is a morphism of Lie algebroids since $\nu$ is holomorphic.

Given a map $\Phi=(\ol{\Phi},\ul{\Phi})$ form $TN$ to $TM$ we have a the bundle map $\Phi'$ from $F$ to $E$ defined by $\ol{\Phi}'(\lambda,w)=(\lambda,\ol{\Phi}(w))$ and $\ul{\Phi}'(\lambda,n)=(\lambda,\ul{\Phi}(n))$. Then $\Phi'$ is a morphism of Lie algebroids if and only if $\ul{\Phi}$ is a holomorphic map and $\ol{\Phi}=T\ul{\Phi}$. Indeed, $\Phi'$ is a morphism if and only if it is admissible $\rho_E\circ\ol{\Phi}'=T\ul{\Phi}'\circ\rho_F$. This is equivalent to
$$
(0,\ol{\Phi}(w)+\lambda J_M(\ol{\Phi}(w)))
=
(0,T\ul{\Phi}(w)+\lambda T\ul{\Phi}(J_N(w)))
$$
for every $(\lambda,w)\in F$, which finally is equivalent to 
$$
\ol{\Phi}=T\ul{\Phi}
\qquand
J_M\circ T\ul{\Phi}=T\ul{\Phi}\circ J_N.
$$

Conversely, if $\Psi$ is a morphism form $F$ to $E$, which is $\lambda$ independent, then $\Psi=\Phi'$ for some holomorphic map $\ul{\Phi}$ and with $\ol{\Phi}=T\ul{\Phi}$. Obviously, $\Phi'$ is a section of $\pi$ if and only is $\ul{\Phi}$ is a section of $\nu$.

\subsection*{Systems with symmetry} The case of a system with symmetry is very important in Physics. We consider a principal bundle $\map{\nu}{P}{M}$ with structure group $G$ and we set $N=M$, $F=TN$ and $E=TP/G$ (the Atiyah algebroid of $P$), with $\pi=([T\nu],\id_M)$. (Here $[T\nu]([v])=T\nu(v)$ for $[v]\in TP/G$.) Sections of $\pi$ are just principal connections on $P$ and a section is a morphism if and only if it is a flat connection. The kernel $K$ is just the adjoint bundle $(P\times\mathfrak{g})/G\to M$. 

An adequate choice of a local basis of sections of $F$, $K$ and $E$ is as follows. Take a coordinate basis $\bar{e_i}=\partial/\partial x^i$ for $F=TM$, take a basis $\{\epsilon_\alpha\}$ of the Lie algebra $\mathfrak{g}$ and the corresponding sections of the adjoint bundle $\{e_\alpha\}$, so that $C^\gamma_{\alpha\beta}$ are the structure constants of the Lie algebra. Finally we take sections $\{e_i\}$ of $E$ such that $e_i$ projects to $\bar{e}_i$, that is, we chose a (local) connection and $e_i$ is the horizontal lift of $\bar{e}_i$. Thus we have $[e_i,e_j]=-\Omega^\alpha_{ij}e_\alpha$ and $[e_i,e_\alpha]=0$.
In this case there are no coordinates $u^A$, and  with the above choice of basis we have that the Euler-Lagrange equations are
\begin{align*}
& \pd{y^\alpha_a}{x^b}-\pd{y^\alpha_b}{x^a}
 +C^\alpha_{\beta\gamma}y^\beta_by^\gamma_a
 =\Omega^\alpha_{ab}\\
&\frac{d}{dx^a}\left(\pd{L}{y^\alpha_a}\right)
-\pd{L}{y^\gamma_a}C^\gamma_{\beta\alpha}y^\beta_a=0.
\end{align*}

In particular, if we chose as connection for the definition of our sections a flat connection (for instance a solution of our variational problem or just a coordinate basis) then this equations reduce to 
\begin{align*}
& \pd{y^\alpha_a}{x^b}-\pd{y^\alpha_b}{x^a}
 +C^\alpha_{\beta\gamma}y^\beta_by^\gamma_a
 =0\\
&\frac{d}{dx^a}\left(\pd{L}{y^\alpha_a}\right)
-\pd{L}{y^\gamma_a}C^\gamma_{\beta\alpha}y^\beta_a=0,
\end{align*}
which are the so called Euler-Poincaré equations~\cite{CaRaSh,CaGaRa}.

The case considered in~\cite{CaRa} where the Lagrangian has only a partial symmetry can be analised in a similar way.


\section{Conclusions and outlook}

We have developed a consistent Field Theory defined over Lie algebroids, finding the Euler-Lagrange equations via a multisymplectic form. We also proved that when the base Lie algebroid is a tangent bundle we have a variational formulation. The admissible variations are determined by the geometry of the problem and are not prescribed in an ad hoc manner. It will be interesting to find a variational formulation for the general case.

Particular cases of our Euler-Lagrange equations are the Euler-Poincaré equations for a system defined on the bundle of connections of a principal bundle, the Lagrange-Poincaré equations for systems defined on a bundle of homogeneous spaces and the Lagrange-Poincaré for systems defined in semidirect products. Our formalism can also describe variational problems defined for holomophic maps.

We have also studied a Hamiltonian formalism and we have proved the equivalence with the Lagrangian formalism in the cases of a hyperregular Lagrangian. Moreover, we have defined a unified formalism which incorporates the relevant aspects of both the Lagrangian and the Hamiltonian formalism.

Among the many advantages of our treatment we mention that it is a multisymplectic theory, that is the field equations are obtained via a multisymplectic equation. This will help in studying the theory of reduction of systems with symmetry and to establish a procedure of reduction by stages. In this respect, the results by \'Sniatycki~\cite{Sniatycki} are relevant.

On the other hand, there are many aspects that has been left out of this paper and can be interesting in it own. For instance, the geometry of our extended notion of jet bundles needs to be studied in more detail, defining the total and contact differential which will allow to write the Euler-Lagrange equations in a simplified way. Systems of partial differential equations are also worth to study and condition for the formal integrability of such systems has to be established. For that, it is obviously necessary to define the concept of higher-order jets.

Finally, in~\cite{MaMeSa} we gave an axiomatic definition of a Lie algebroid structure over an affine bundle which encodes the geometric structure necessary for developing time-dependent Classical Mechanics. It would be nice to isolate the geometric structure necessary for developing field theories defining what could be called a Lie multialgebroid.

\end{document}